\documentclass{amsart}

\usepackage{amscd}
\usepackage{amssymb}
\usepackage{amsmath}
\usepackage{amsfonts}
\usepackage{bm}
\usepackage{enumerate}
\usepackage[all,cmtip]{xy}
\usepackage{varioref}
\usepackage{tikz-cd}
\usepackage{setspace}
\usepackage[colorlinks,linktocpage]{hyperref}
\hypersetup{citecolor=blue,linkcolor=blue}

\newtheorem{thm}{Theorem}[section]
\newtheorem{prop}[thm]{Proposition}
\newtheorem{lem}[thm]{Lemma}
\newtheorem{cor}[thm]{Corollary}
\newtheorem{conj}[thm]{Conjecture}

\theoremstyle{definition}

\theoremstyle{remark}
\newtheorem{rem}[thm]{Remark}
\newtheorem*{nota}{Notation}

\newcommand{\kk}{\Bbbk}

\newcommand{\ZZ}{\mathcal{Z}}

\newcommand{\A}{\mathcal{A}}
\newcommand{\GG}{\mathcal{G}}
\newcommand{\JJ}{\mathcal{J}}

\newcommand{\PP}{\mathcal{P}}
\newcommand{\QQ}{\mathcal{Q}}
\newcommand{\RR}{\mathcal{R}}
\newcommand{\MM}{\mathcal{M}}
\newcommand{\II}{\mathcal{I}}
\newcommand{\DD}{\mathcal{D}}

\newcommand{\Zz}{\mathbb{Z}}
\newcommand{\Ff}{\mathbb{F}}
\newcommand{\Cc}{\mathbb{C}}
\newcommand{\Rr}{\mathbb{R}}

\newcommand{\lk}{\mathrm{lk}}
\newcommand{\st}{\mathrm{star}}
\newcommand{\xr}{\xrightarrow}
\newcommand{\wh}{\widehat}
\newcommand{\w}{\widetilde}

\DeclareMathOperator{\im}{Im}

\DeclareMathOperator{\row}{row}
\DeclareMathOperator{\Hom}{Hom}
\DeclareMathOperator{\Ext}{Ext}
\DeclareMathOperator{\supp}{supp}
\DeclareMathOperator*{\colim}{colim}

\numberwithin{equation}{section}

\begin{document}
\title[On the integral cohomology of real toric manifolds]{On the integral cohomology of real toric manifolds}
\author[F.~Fan]{Feifei Fan}
\thanks{This work was supported by the National Natural Science Foundation of China (Grant No. 12271183) and by the Guangdong Basic and Applied Basic Research Foundation (Grant No. 2023A1515012217).}
\address{Feifei Fan, School of Mathematical Sciences, South China Normal University, Guangzhou, 510631, China.}
\email{fanfeifei@mail.nankai.edu.cn}

\begin{abstract}
	Real toric manifolds are the real loci of nonsingular complete toric varieties. We compute the integral cohomology groups of real toric manifolds from the combinatorial data encoded by the underlying simplicial fans, generalizing a formula due to Cai and Choi. Furthermore, we determine the image of the integral cohomology in the $\Zz/2$-cohomology under the mod $2$ reduction map. As an application, a combinatorial-algebraic criterion is established for when a real toric manifold is $\mathrm{spin}^c$. 
	We also describe the product structure of their integral cohomology rings. In particular, we establish a simple combinatorial formula for the quotient of the cohomology ring by the ideal consisting of elements annihilated by $2$. 
\end{abstract}

\maketitle

\section{Introduction}\label{sec:introduction}
Toric varieties in algebraic geometry have deep connections with  combinatorics, commutative algebra, symplectic geometry, and topology. A \emph{toric variety} $X$ of dimension $n$ is a normal complex algebraic variety containing an algebraic torus $(\Cc^*)^n$ as a Zariski open subset, such that the action of the torus on itself extends to the whole variety. 
Nonsingular complete (smooth and compact in the usual topology) toric varieties are also called \emph{toric manifolds}. The classification theorem of toric varieties says that there is a bijective correspondence between the isomorphism classes of complex
$n$-dimensional toric manifolds and complete nonsingular fans in $\Rr^n$. In particular, a projective toric manifold
$X$ is defined by the normal fan of a \emph{lattice Delzant $n$-polytope} $P_X\subset \Rr^n$.

In \cite{DJ91}, Davis and Januszkiewicz introduced topological generalizations of projective toric manifolds, now known as \emph{quasitoric manifolds}. A quasitoric manifold is a closed smooth $2n$-manifold $M$ with a locally standard $T^n$-action (i.e., it is locally isomorphic to the standard action of $T^n=(S^1)^n$ on $\Rr^{2n}$) such that the orbit space $M/T^n$ can be identified with a simple $n$-polytope. A toric manifold is not necessarily a quasitoric
manifold and vice versa. However, both toric and quasitoric manifolds are examples of \emph{topological toric manifolds} introduced in \cite{IFM13}. A topological toric manifold  is a closed smooth $2n$-manifold with an effective smooth
action of $(\Cc^*)^n$ having an open dense orbit and covered by finitely many invariant open subsets each equivariantly diffeomorphic to a direct sum of complex one-dimensional smooth representation spaces of $(\Cc^*)^n$.

Topological toric manifolds have several favorable topological properties. For example, their singular cohomology can be described nicely as a quotient ring of a polynomial ring. The formula for the integral cohomology rings of toric manifolds was obtained by Jurkiewicz \cite{Jur80} (for the projective case) and Danilov \cite{Dan78} (for the general case), and similar formulas were obtained for quasitoric manifolds in \cite{DJ91} and for topological toric manifolds in \cite{IFM13}. 

In toric topology, there are also real analogs of these toric spaces. For example, the real locus $X^{\Rr}$ of a toric manifold $X$ of complex dimension $n$ is a nonsingular real variety of dimension $n$, called a \emph{real toric manifold}. The real analog of a quasitoric manifold is a closed smooth manifold of dimension $n$ with a
locally standard action of $(\Zz/2)^n=(S^0)^n\subset (\Rr^*)^n$, called a \emph{small cover} (\cite{DJ91}). Also, a closed smooth $n$-manifold with an
effective smooth action of $(\Rr^*)^n$ is called a \emph{real topological toric manifold} if
one of its orbits is open and dense and it is covered by finitely many invariant
open subsets each equivariantly diffeomorphic to a direct sum of
real one-dimensional smooth representation spaces of $(\Rr^*)^n$  (\cite{IFM13}).
As in the complex case, both real toric manifolds and small covers are examples of  
real topological toric manifolds.

More generally, let $K$ be an abstract simplicial complex on the vertex set $[m]=\{1,2,\dots,m\}$. The \emph{real moment-angle complex} $\Rr\ZZ_K$ associated to $K$ is the topological space defined as 
\begin{equation}\label{eq:m-a}
	\Rr\ZZ_K=\bigcup_{\sigma\in K}\{(x_1,\dots,x_m)\in (D^1)^m\mid x_i\in S^0\text{ when }i\not\in\sigma\},
\end{equation}
where $D^1$ is the interval $[-1,1]$  and $S^0=\{-1,1\}$ is its boundary.
It should be mentioned that $\Rr\ZZ_K$ is a PL (piecewise linear) manifold if and only if $K$ is a PL sphere (see \cite[Theorem 2.3]{C17}). Furthermore, if $K$ is a \emph{star-shaped} simplicial sphere (a triangulated sphere isomorphic to the underlying complex of a complete simplicial fan), then $\Rr\ZZ_K$ is a smooth manifold (see \cite[Theorem 6.5.2]{BP15}). 
There is a canonical $(\Zz/2)^m$-action on $\Rr\ZZ_K$ which comes from the $\Zz/2$-action on the interval $[-1,1]$ by reflection. Let $G\subset (\Zz/2)^m$ be a subgroup acting freely on $\Rr\ZZ_K$. Then the quotient space $\Rr\ZZ_K/G$ is called a \emph{real toric space}. Any subgroup of $(\Zz/2)^m$
can be realized as the kernel of a surjective linear map $\Lambda:(\Zz/2)^m\to (\Zz/2)^n$. One can show (see \cite[Lemma 3.1]{CKT17}) that $\ker\Lambda$ acts freely on $\Rr\ZZ_K$ if and only if $\Lambda$ satisfies the \emph{nonsingularity condition}:
\[\Lambda(\boldsymbol{e}_{i_1}),\dots,\Lambda(\boldsymbol{e}_{i_k})\text{ are linearly independent in }(\Zz/2)^n\text{ if }\{i_1,\dots,i_k\}\in K,\]
where $\boldsymbol{e}_i$ is the $i$th basis vector of $(\Zz/2)^m$. In this case, we denote the real toric space $\Rr\ZZ_K/\ker\Lambda$ by $M(K,\Lambda)$. 

It is shown in \cite{IFM13} that a real topological toric manifold is a real toric space $M(K,\Lambda)$ for some star-shaped simplicial sphere $K$. In particular, an $n$-dimensional real toric manifold $X^{\Rr}$ is  homeomorphic to the real toric space $M(K_X,\Lambda_X)$, where $K_X$ is the star-shaped $(n-1)$-sphere corresponding to the underlying fan of $X$ and $\Lambda_X:(\Zz/2)^m\to (\Zz/2)^n$ is the map given by the  primitive generators (mod $2$) of the $m$ rays of the underlying fan (see \cite[Ch. 5]{BP15}).

Davis and Januszkiewicz \cite{DJ91} have established the formula for the $\Zz/2$-cohomology ring of the real toric space $M(K,\Lambda)$ 
 when $K$ is a simplicial sphere or more generally a Cohen-Macaulay complex over $\Zz/2$ and $\Lambda$ has the smallest possible rank, i.e., $\dim K+1$ (see Theorem \ref{thm:DJ}). Such real toric spaces include all real topological toric manifolds. However, relatively little is known about the integral cohomology rings of these real toric spaces. 
Only the additive structure of the integral cohomology of small covers and real toric manifolds associated to shellable fans was obtained by Cai and Choi in \cite{CC21}. (The rational Betti numbers of small covers were computed by Suciu and Trevisan in an earlier work \cite{Suciu13}, \cite[Sec. VI.2]{Trev12}.)

In \cite{CP17,CP20}, Choi and Park gave an explicit formula for the $\kk$-cohomology rings of general real toric spaces under the assumption that $2$ is invertible in the coefficient ring $\kk$, where the multiplicative structure has an especially beautiful combinatorial description in terms of $K$ and $\Lambda$. Later, for any real toric space $M(K,\Lambda)$ and any coefficient ring $\kk$, Franz \cite{Franz22} constructed a DGA (differential graded algebra) over $\kk$, whose cohomology is isomorphic to $H^*(M(K,\Lambda);\kk)$ as $\kk$-algebras. Each of these results has its own advantages and disadvantages. For example, Franz's result has the advantage of applying to arbitrary coefficient rings, but the cohomology of the DGA he constructed is itself difficult to compute in general.  

In this paper, we focus on the integral cohomology of real toric spaces $M(K,\Lambda)$, especially when $K$ is a PL sphere and $\Lambda$ has rank $\dim K+1$, so that $M(K,\Lambda)$ is a manifold of dimension $\dim K+1$. As explained above, any real (topological) toric manifold is such a manifold since star-shaped spheres are PL.  
 
First, for any real toric space $M(K,\Lambda)$, we construct a free differential graded $\Zz$-module $(\RR(K,\Lambda),d)$ (coming from the canonical cellular decomposition of $\Rr\ZZ_K$), such that $H^*(\RR(K,\Lambda),d)$ is isomorphic to $H^*(M(K,\Lambda);\Zz)$ as groups (Theorem \ref{thm:isomorphism}). 
We also construct a DGA model (Theorem \ref{thm:DGA}) for real toric spaces that differs from Franz’s model \cite{Franz22}. Then, we extend the result of Cai and Choi \cite{CC21} by determining the additive structure of $H^*(M(K,\Lambda);\Zz)$ when $K$ is a PL sphere and $\Lambda$ has rank $\dim K+1$ (Proposition \ref{prop:group}). In this situation, the image of the mod $2$ reduction map $H^*(M(K,\Lambda);\Zz)\to H^*(M(K,\Lambda);\Zz/2)$ can be determined (Theorem \ref{thm:mod 2 image}). As an application, we establish a criterion for a real topological toric manifold to admit a $\mathrm{spin}^c$ structure (Proposition \ref{prop:spinc}). Furthermore, with these assumptions on $K$ and $\Lambda$ we can describe the multiplicative structure of the integral cohomology ring $H^*(M(K,\Lambda);\Zz)$ (Theorem \ref{thm:ring isomorphism}). In particular, we establish a combinatorial formula for the quotient ring of $H^*(M(K,\Lambda);\Zz)$ modulo the ideal consisting of elements annihilated by $2$ (Theorem \ref{thm:ring}), which is analogous to the formula of Choi and Park \cite[Main Theorem]{CP20} for $\kk$-cohomology with $1/2\in\kk$. 

\begin{nota}
	For a given space $X$, we write $H^*(X)$ ($\w H^*(X)$) for the (reduced) integral cohomology $H^*(X;\Zz)$. Given two $\Zz$-modules $A,B$, we write $\Hom_{\Zz}(A,B)$ simply as $\Hom(A,B)$ and write $A\otimes_\Zz B$ as $A\otimes B$.
\end{nota}

\section{Main results}\label{sec:result}
Throughout this paper, we assume that $K$ is a simplicial complex on the set $[m]=\{1,2,\dots,m\}$, possibly with ghost vertices. For a non-empty face $\sigma$ of $K$, we set $\partial\sigma=2^\sigma\setminus\{\sigma\}$ to be the boundary complex of the simplex $2^\sigma=\{\tau:\tau\subset\sigma\}\subset K$. For a face $\sigma\in K$, the \emph{link} and the \emph{star} of $\sigma$ in $K$ are the subcomplexes
\begin{gather*}
	\lk_K\sigma=\{\tau\in K\mid\tau\cup\sigma\in K,\tau\cap\sigma=\emptyset\},\\
	\st_K\sigma=\{\tau\in K\mid\tau\cup\sigma\in K\}.
\end{gather*}
The \emph{join} of two simplicial complexes $K$ and $K'$ with  disjoint vertex sets  is the simplicial complex
\[K*K'=\{\sigma\cup\sigma'\mid\sigma\in K,\,\sigma'\in K'\}.\]
The \emph{stellar subdivision} of $K$ at a non-empty face $\sigma$ is the simplicial complex
\[S_\sigma K=(K\setminus\st_K\sigma)\cup (\{v_\sigma\}*\partial\sigma*\lk_K\sigma),\]
where $v_\sigma$ is a vertex not in $K$.

Fix a field $\Ff$ and let $\Ff[x_1,\dots,x_m]$ be the graded polynomial ring with $\deg x_i=1$. For a subset $\sigma\subset [m]$, denote by $x_\sigma$ the square-free monomial $\prod_{i\in\sigma}x_i$. The \emph{Stanley-Reisner ring} $\Ff[K]$ of the simplicial complex $K$ is the quotient of $\Ff[x_1,\dots,x_m]$ by a square-free monomial ideal
\[\Ff[K]=\Ff[x_1,\dots,x_m]/\langle x_\sigma\mid\sigma\not\in K\rangle.\]
A simplicial complex $K$ is \emph{Cohen-Macaulay over $\Ff$} if $\Ff[K]$ is a \emph{Cohen-Macaulay ring}, i.e., the maximum length of a regular
sequence in $\Ff[K]$ equals $\dim K+1$. By a classical result of Reisner \cite{Rei76}, $\Ff[K]$ is Cohen-Macaulay if and only if $\w H^i(\lk_K\sigma;\Ff)=0$ for all $\sigma\in K$ (including $\sigma=\emptyset$) and $i<\dim\lk_K\sigma$. For instance, simplicial spheres are Cohen-Macaulay over any field $\Ff$.

Let $(C_*(K),\partial)$ ($(\w C_*(K),\partial)$) denote the (augmented) simplicial chain complex of $K$, and let $(C^*(K),\delta)$ ($(\w C^*(K),\delta)$) denote its dual cochain complex. For an oriented face $\sigma\in K$, whose orientation is induced by an ordering of its vertices, denote by $\sigma^*\in\w C^*(K)$ the basis cochain dual to $\sigma$; it takes value $1$ on $\sigma$ and vanishes on all other faces. For a subset $\omega\subset [m]$, denote by $K_\omega=\{\sigma\in K\mid \sigma\subset \omega\}$ the full subcomplex of $K$ with respect to $\omega$.
Note that there is a bijection between $(\Zz/2)^m$ and the power set $2^{[m]}$: \[(b_1,\dots,b_m)\mapsto \{i\in[m]\mid b_i\neq 0\}.\]
So we do not distinguish an element of $(\Zz/2)^m$ and the corresponding subset of $[m]$. 

Following the construction in \cite{C17}, let $\Zz\langle u_1,\dots,u_m;t_1,\dots,t_m\rangle$ be the differential graded free $\Zz$-algebra with $2m$ generators such that 
\[\deg u_i=1,\ \deg t_i=0;\ dt_i=u_i,\ du_i=0,\]
and let $\RR$ be the quotient of $\Zz\langle u_1,\dots,u_m;t_1,\dots,t_m\rangle$ under the following relations:
\begin{gather*}
u_it_i=u_i,\ \ t_iu_i=0,\ \ u_it_j=t_ju_i,\ \ t_it_i=t_i,\\
u_iu_i=0,\ \ u_iu_j=-u_ju_i,\ \ t_it_j=t_jt_i
\end{gather*}
for $i,j\in[m]$ and $i\neq j$. For an ordered subset $\sigma=\{i_1,\dots,i_k\}\subset[m]$, denote by $u_\sigma$ (resp. $t_\sigma$) the monomial $u_{i_1}\cdots u_{i_k}$ (resp. $t_{i_1}\cdots t_{i_k}$). For a simplicial complex $K$, let $\II_K$ be the ideal of $\RR$ generated by all square-free monomials $u_\sigma$ such that $\sigma$ is not a face of $K$. We write  $\RR_K=\RR/\II_K$. 
It is easy to see that as a $\Zz$-module, $\RR_K$ is freely generated by the square-free monomials $u_\sigma t_{\omega\setminus\sigma}$, where $\sigma\subset\omega\subset [m]$ and $\sigma \in K$. As we will see in Section \ref{sec:m-a}, $\RR_K$ is actually a DGA model for the cellular cochain complex of the real moment-angle complex $\Rr\ZZ_K$.

For each $g=(g_1,\dots,g_m)\in(\Zz/2)^m$, define a $\Zz$-module endomorphism  
\begin{equation}\label{eq:action}
	\begin{split}
		\phi_g:\RR&\to \RR\\
		u_\sigma t_{\omega\setminus\sigma}&\mapsto (-1)^{|g\cap\sigma|}u_\sigma t_{\omega\setminus\sigma}',
	\end{split}
\end{equation}	
where $\sigma\subset\omega\subset [m]$, $\sigma \in K$, and
\[t_{\omega\setminus\sigma}'=\prod_{i\in \omega\setminus\sigma}t_i'\ \text{ and }\ t_i'=\begin{cases}
	t_i&\text{ if }g_i=0,\\
	1-t_i&\text{ if }g_i=1.
\end{cases}
\]
Then the group $(\Zz/2)^m$ acts on $\RR$ by $g\cdot\alpha=\phi_g(\alpha)$, $g\in (\Zz/2)^m$, and this action passes to $\RR_K$. It should be noted that this action of $(\Zz/2)^m$ does not preserve the product structure. For instance, if $i\in g$ then $\phi_g(t_iu_i)=\phi_g(0)=0$, while $\phi_g(t_i)\phi_g(u_i)=(1-t_i)(-u_i)=-u_i$.

Let $\Lambda:(\Zz/2)^m\to (\Zz/2)^n$ be a surjective linear map. Fixing a basis of $(\Zz/2)^n$, the map $\Lambda$ can be represented by an $n\times m$ $\Zz/2$-matrix $(\lambda_{ij})$, also denoted $\Lambda$. (Changing basis of $(\Zz/2)^n$ changes $\Lambda$ to $M\Lambda$ for some $M\in GL_n(\Zz/2)$.) We denote by $\row\Lambda$ the subspace of $(\Zz/2)^m$ spanned by the row
vectors of the matrix $\Lambda$. In particular, if $\Lambda$ satisfies the 
nonsingularity condition for $K$ introduced in Section \ref{sec:introduction}, then $\Lambda$ is called a (mod 2) \emph{characteristic function} over $K$. As we mentioned earlier, $\ker\Lambda$ acts freely on $\Rr\ZZ_K$ if and only if $\Lambda$ is a characteristic function  over $K$, and the quotient space $M(K,\Lambda)=\Rr\ZZ_K/\ker\Lambda$
is called a \emph{real toric space} in this case. 

For a surjective linear map $\Lambda:(\Zz/2)^m\to (\Zz/2)^n$, let $\RR(K,\Lambda)=(\RR_K)^{\ker\Lambda}$ be the submodule of $\RR_K$ consisting of  elements fixed by the subgroup $\ker\Lambda\subset(\Zz/2)^m$.
Then $\RR(K,\Lambda)$ is a differential graded free $\Zz$-module. Write $\RR^i(K,\Lambda)$ for the graded piece of $\RR(K,\Lambda)$ in degree $i$. Notice that $\RR(K,\Lambda)$ is not a DGA, since it is not closed under multiplication. 

\begin{thm}\label{thm:isomorphism}
	Let $\Lambda:(\Zz/2)^m\to (\Zz/2)^n$ be a characteristic function over $K$. Then  $\RR(K,\Lambda)\subset\RR_K$ is the submodule
	\[\{\sum_{g\in\ker\Lambda}g\cdot\alpha\mid \alpha\in\RR_K\},\]
	and there is an isomorphism of graded $\kk$-modules for any coefficient ring $\kk$
	\[H^*(M(K,\Lambda);\kk)\cong H^*(\RR(K,\Lambda)\otimes\kk).\]
\end{thm}

The isomorphism in Theorem \ref{thm:isomorphism} is not multiplicative, since $\RR(K,\Lambda)$ is not a DGA. Nevertheless, we can construct a DGA model for real toric spaces by using a different DGA model for $\Rr\ZZ_K$. 

We first define a $\Zz$-DGA $\Gamma$ as follows. As a $\Zz$-module, $\Gamma$ is freely generated by the $(n+1)$-tuples $[g_0,\dots,g_n]$, $n\geq 0$, $g_i\in\Zz/2$, with $g_i\neq g_{i+1}$ for $0\leq i\leq n-1$, and $\deg[g_0,\dots,g_n]=n$. The product structure is given by 
\[[g_0,\dots,g_k]\cdot[g_0',\dots,g_l']=\begin{cases}
[g_0,\dots,g_k,g_1',\dots,g_{l}'],&\text{if }g_k=g_0',\\
0,&\text{otherwise,}
\end{cases}\]
and the differential is given by
\[d([g_0,\dots,g_n])=[g_{-1},g_0,\dots,g_n]+(-1)^{n+1}[g_0,\dots,g_n,g_{n+1}],\]
where $g_{-1}$ (resp. $g_{n+1}$) is the unique element of $\Zz/2\setminus\{g_0\}$ (resp. $\Zz/2\setminus\{g_n\}$). It is straightforward to verify that the differential $d$ satisfies the Leibniz rule. ($\Gamma$ is essentially the dual of the chain complex $C(E\Zz/2)$ constructed in \cite[Sec. 2.3]{Franz22}.)

Since $\Zz/2=\{1,-1\}$ (here we use the multiplicative representation of $\Zz/2$), $\Gamma$ has only two additive generators in each degree $n$:
\[\mu_n=[-1,1,\dots,(-1)^{n+1}]\ \text{ and }\ \mu_n'=[1,-1,\dots,(-1)^n].\]
With this notation, one easily checks that $\mu_0+\mu_0'$ is the identity element of $\Gamma$ and there are the following relations:
\begin{gather}\label{eq:relations}
	\begin{gathered}
	\mu_i\mu_j=\begin{cases}
		\mu_{i+j},&i\text{ even},\\
		0,&i\text{ odd},
	\end{cases}\quad 
	\mu'_i\mu'_j=\begin{cases}
	\mu'_{i+j},&i\text{ even},\\
	0,&i\text{ odd},
	\end{cases}\\
	\mu_i\mu_j'=\begin{cases}
		0,&i\text{ even},\\
		\mu_{i+j},&i\text{ odd},
	\end{cases}\quad
	\mu_i'\mu_j=\begin{cases}
		0,&i\text{ even},\\
		\mu_{i+j}',&i\text{ odd},
	\end{cases}\\
	d(\mu_i)=\mu_{i+1}'-(-1)^i\mu_{i+1},\quad d(\mu_i')=\mu_{i+1}-(-1)^i\mu_{i+1}'.
	\end{gathered}
\end{gather}

The algebra $\Gamma$ admits a $\Zz/2$-action given by
\[g\cdot[g_0,\dots,g_n]=[gg_0,\dots,gg_n],\ \ g\in\Zz/2.\] 
This action is clearly algebraic.

Let $\QQ$ be the $m$-fold tensor product $\Gamma^{\otimes m}$ with the standard product formula
\[(\alpha_1\otimes\cdots\otimes\alpha_m)(\beta_1\otimes\cdots\otimes\beta_m)=(-1)^s\alpha_1\beta_1\otimes\cdots\otimes\alpha_m\beta_m,\]
where $s=\sum_{i=1}^m\deg\beta_i\sum_{j>i}\deg\alpha_j$, and with the standard differential formula
\[d(\alpha_1\otimes\cdots\otimes\alpha_m)=\sum_{i=1}^m(-1)^{\sum_{j=1}^{i-1}\deg \alpha_j}\alpha_1\otimes\cdots \otimes d(\alpha_i)\otimes\cdots\otimes\alpha_m.\]
Then $\QQ$ is a DGA, and it has a refined $\Zz^m$-grading. For each $i\in[m]$, define 
\begin{equation}\label{eq:notation}
\begin{gathered}
	x_{n,i}=1\otimes\cdots\otimes1\otimes \mu_n\otimes1\otimes\cdots\otimes1,\\ y_{n,i}=1\otimes\cdots\otimes1\otimes \mu_n'\otimes1\otimes\cdots\otimes1,
\end{gathered}
\end{equation}
where $\mu_n$ or $\mu_n'$ occupies the $i$th entry.  Given an ordered subset $\sigma=\{i_1,\dots,i_k\}\subset[m]$, set
\[y_\sigma=\prod_{j=1}^k y_{0,i_j},\quad z_\sigma=\prod_{j=1}^k (x_{1,i_j}-y_{1,i_j}).\]

For a homogeneous basis element $\alpha=\alpha_1\otimes\cdots\otimes\alpha_m\in\QQ$ with respect to the $\Zz^m$-grading, define the support of $\alpha$
\[\supp\alpha=\{i\in[m]\mid\deg\alpha_i>0\}.\] 
Let $\JJ_K$ be the ideal of $\QQ$ generated by all homogeneous basis elements  $\alpha$ such that $\supp\alpha$ is not a face of $K$, and let $\QQ_K=\QQ/\JJ_K$. It can be shown that $\QQ_K$ is also a DGA model for $\Rr\ZZ_K$ (Propositions \ref{prop:DGA model} and \ref{cor:E-Z}). 

The algebra $\QQ$ admits a $(\Zz/2)^m$-action given by
\[g\cdot(\alpha_1\otimes\cdots\otimes\alpha_m)=(g_1\cdot\alpha_1)\otimes\cdots\otimes (g_m\cdot\alpha_m),\ \ g=(g_1,\dots,g_m),\]
which passes to $\QQ_K$. Since this action preserves the DGA structure, the $G$-invariants $(\QQ_K)^G$ is a DGA for any subgroup $G\subset(\Zz/2)^m$.

\begin{thm}\label{thm:DGA}
Let $\Lambda:(\Zz/2)^m\to (\Zz/2)^n$ be a characteristic function over $K$, and let $\QQ(K,\Lambda)=(\QQ_K)^{\ker\Lambda}$. Then $\QQ(K,\Lambda)\subset\QQ_K$ is the subring 
\[\{\sum_{g\in\ker\Lambda}g\cdot\alpha\mid \alpha\in\QQ_K\},\]
and there is an isomorphism of graded $\kk$-algebras for any coefficient ring $\kk$
\[H^*(M(K,\Lambda);\kk)\cong H^*(\QQ(K,\Lambda)\otimes\kk).\]	
\end{thm}

\begin{rem}
The explicit computation of $H^*(\RR(K,\Lambda)\otimes\kk)$ or $H^*(\QQ(K,\Lambda)\otimes\kk)$ is difficult for general $(K,\Lambda)$ and $\kk$ except when $M(K,\Lambda)=\Rr\ZZ_K$ so that $\RR(K,\Lambda)=\RR_K$ or when $1/2\in\kk$. For these special cases, the $\kk$-module structure of $H^*(M(K,\Lambda);\kk)$ (and even the ring structure of $H^*(M(K,\Lambda);\kk)$) has a nice combinatorial description. (See \cite[Proposition 3.3]{C17}, \cite[Proposition A.2]{FMW24}, \cite[Corollary 4.10]{Yu24} for $M(K,\Lambda)=\Rr\ZZ_K$ and \cite[Theorem 4.6]{CP17}, \cite[Main Theorem]{CP20} for $1/2\in\kk$.)
\end{rem}

For any $\omega\in\row\Lambda$ and any oriented face $\sigma\in K_\omega$, define
\begin{gather}\label{eq:f,h}
	\begin{gathered}
		f_{\sigma,\omega}=\sum_{\omega'\subset\omega\setminus\sigma}(-2)^{|\omega'|}u_\sigma t_{\omega'}\in \RR_K,\\
		h_{\sigma,\omega}=\sum_{\omega'\subset\omega\setminus\sigma}(-2)^{|\omega'|}z_\sigma y_{\omega'}\in \QQ_K.
	\end{gathered}
\end{gather}
It can be shown that (Lemma \ref{lem:invariants}) $f_{\sigma,\omega}\in \RR(K,\Lambda)$, $h_{\sigma,\omega}\in \QQ(K,\Lambda)$ and that there are cochain maps:
\begin{gather}\label{eq:map phi_omega}
		\begin{gathered}
	\varphi_\omega:(\w C^*(K_\omega),2\delta)\to (\RR(K,\Lambda),d),\ \sigma^*\mapsto (-1)^{|\sigma|}f_{\sigma,\omega},\\
	\psi_\omega:(\w C^*(K_\omega),2\delta)\to (\QQ(K,\Lambda),d),\ \sigma^*\mapsto (-1)^{|\sigma|}h_{\sigma,\omega}.
		\end{gathered}
\end{gather}
Notice that we use $2\delta$ rather than $\delta$ as the differential operator of the cochain complex on the left side. Let 
\[\A^*(K_\omega)=(\w C^{*-1}(K_\omega),2\delta)\ \text{ and }\  \A(K,\Lambda)=\bigoplus_{\omega\in \row\Lambda}\A^*(K_\omega).\] 
Then there are induced homomorphisms of cohomology groups: 
\begin{gather}\label{eq:map}
	\begin{gathered}
	\Phi_{\Lambda}:H^*(\A(K,\Lambda))\cong \bigoplus_{\omega\in \row\Lambda} H^*(\A^*(K_\omega))\to H^*(\RR(K,\Lambda)),\\
	\Psi_{\Lambda}:H^*(\A(K,\Lambda))\cong \bigoplus_{\omega\in \row\Lambda} H^*(\A^*(K_\omega))\to H^*(\QQ(K,\Lambda)),
	\end{gathered}
\end{gather}
where the restriction of $\Phi_{\Lambda}$ (resp. $\Psi_{\Lambda}$) to $H^*(\A^*(K_\omega))$ is induced by $\varphi_\omega$  (resp. $\psi_\omega$).

One easily checks that the $\Zz$-module homomorphism 
\[\RR_K\to\QQ_K,\ \ u_\sigma t_{\omega\setminus\sigma}\mapsto z_\sigma y_{\omega\setminus\sigma},\] 
where $\sigma\in K$ and $\sigma\subset\omega\subset[m]$, is a $(\Zz/2)^m$-equivariant cochain map. Hence, there is a commutative diagram:
\begin{equation}\label{diag:Phi and Psi}
	\begin{gathered}
		\xymatrix{
			H^*(\A(K,\Lambda))\ar[r]^-{\Phi_\Lambda}\ar@{=}[d]&H^*(\RR(K,\Lambda))\ar[d]\\
			H^*(\A(K,\Lambda))\ar[r]^-{\Psi_\Lambda}&H^*(\QQ(K,\Lambda))}
	\end{gathered}
\end{equation}

When $K$ is a Cohen-Macaulay complex over $\Zz/2$ and $\Lambda$ has rank $\dim K+1$ we can describe the image of $H^*(\A(K,\Lambda))$ in the mod $2$ cohomology of $M(K,\Lambda)$, while the latter is determined by Davis and Januszkiewicz as follows.

\begin{thm}[{\cite[Theorem 5.12]{DJ91}}]\label{thm:DJ}
	Let $K$ be an $(n-1)$-dimensional Cohen-Macaulay complex over $\Zz/2$, and assume that $\Lambda:(\Zz/2)^m\to (\Zz/2)^n$ is a characteristic function over $K$. Let $\lambda_i\in\row\Lambda$ be the $i$th row vector and let $\Theta_\Lambda$ be the ideal of $\Zz/2[K]$ generated by the linear forms
	\[\theta_i=\sum_{j\in\lambda_i}x_j,\ \ 1\leq i\leq n.\]
	Then the mod $2$ cohomology ring of $M(K,\Lambda)$ is given by
	\[H^*(M(K,\Lambda);\Zz/2)\cong \Zz/2[K]/\Theta_\Lambda.\]
\end{thm}

\begin{thm}\label{thm:mod 2 image}
	Let $K$ and $\Lambda$ be as in Theorem \ref{thm:DJ}. Then for a cohomology class $\alpha=[\sum_{\sigma\in K_\omega}k_\sigma\sigma^*]\in H^*(\A^*(K_\omega))$ with $\omega\in \row\Lambda$, the mod $2$ reduction of $\Phi_\Lambda(\alpha)$ (or $\Psi_\Lambda(\alpha)$) in $H^*(M(K,\Lambda);\Zz/2)$ is the polynomial
	\[P_\alpha:=\sum_{\sigma\in U}x_\sigma,\ \text{ where }\ U=\{\sigma\in K_\omega\mid k_\sigma \text{ is odd}\}.\]
\end{thm}

When $K$ is a PL sphere and $\Lambda$ has rank $\dim K+1$, the integral cohomology group of $M(K,\Lambda)$ has a nice description.
\begin{thm}\label{thm:group}
Let	$K$ be an $(n-1)$-dimensional PL sphere, and assume that $\Lambda:(\Zz/2)^m\to (\Zz/2)^n$ is a characteristic function over $K$. Then the maps $\Phi_\Lambda$, $\Psi_\Lambda$ in \eqref{eq:map}
are surjective. Moreover, let $\A_2$, $\RR_2$, $\QQ_2$ be the subgroups of $H^*(\A(K,\Lambda))$, $H^*(\RR(K,\Lambda))$, $H^*(\QQ(K,\Lambda))$, respectively, consisting of elements annihilated by $2$. Then there are graded isomorphisms of quotient groups
\begin{gather*}
	\Phi_\Lambda:H^*(\A(K,\Lambda))/\A_2\cong H^*(\RR(K,\Lambda))/\RR_2,\\
	\Psi_\Lambda:H^*(\A(K,\Lambda))/\A_2\cong H^*(\QQ(K,\Lambda))/\QQ_2.
	\end{gather*}
\end{thm}
Combining Theorems \ref{thm:isomorphism} and \ref{thm:group} gives the following corollary which extends the result of Cai and Choi \cite[Corollary 5.3]{CC21} for shellable spheres to PL spheres. It is known that the inclusion $\{\text{shellable spheres}\}\subset \{\text{PL spheres}\}$ is strict in every dimension greater than 2 (see \cite{Lic91}), so this is a substantial extension. 
\begin{cor}\label{cor:group}
With the same hypotheses as Theorem \ref{thm:group}, let $\GG^i(K,\Lambda)$ be the group $\bigoplus_{\omega\in \row\Lambda}\w H^{i-1}(K_\omega)$, $i\geq 0$ (here we use the convention $\w H^{-1}(\emptyset)=\Zz$), and let $p$ be an odd prime. Then
\begin{enumerate}[(a)]
	\itemsep=3pt
	\item\label{item:1} The $\Zz$-summands (resp. $\Zz/{p^k}$-summands, $k\geq 1$) of $\GG^i(K,\Lambda)$ are in one-to-one correspondence with those of $H^i(M(K,\Lambda))$.
	\item\label{item:2}	For $k\geq 1$, the $\Zz/{2^k}$-summands of $\GG^i(K,\Lambda)$ are in one-to-one correspondence with the $\Zz/{2^{k+1}}$-summands of $H^i(M(K,\Lambda))$.
\end{enumerate}
\end{cor}
\begin{proof}
Part \eqref{item:1} follows directly from Theorems \ref{thm:isomorphism} and  \ref{thm:group}. For part \eqref{item:2}, first note that the $\Zz/{2^k}$-summands of $\GG^i(K,\Lambda)$ are in one-to-one correspondence with the $\Zz/{2^{k+1}}$-summands of $H^i(\A(K,\Lambda))$ since the differential of $\A*(K_\omega)$ is twice the differential of $\w C^*(K_\omega)$. By the last assertion of Theorem \ref{thm:group} for $\Phi_\Lambda$ and Lemma \ref{lem:truncated}, $2H^i(\A(K,\Lambda))\cong 2H^i(\RR(K,\Lambda))$. Hence, for $k\geq 1$, the $\Zz/{2^{k+1}}$-summands of $H^i(\A(K,\Lambda))$ must be in one-to-one correspondence with those of $H^i(\RR(K,\Lambda))$.
 Then part \eqref{item:2} follows by Theorem \ref{thm:isomorphism}.
\end{proof}
By a result of Stanley (see \cite[Ch. II]{S96}), if $K$ is an $(n-1)$-dimensional \emph{Cohen-Macaulay complex} over $\Zz/2$ and $\Lambda:(\Zz/2)^m\to (\Zz/2)^n$ is a characteristic function over $K$, then the $i$th $\Zz/2$-Betti number of $M(K,\Lambda)$ equals the $i$th component of the \emph{$h$-vector} of $K$, a sequence of combinatorial numbers of $K$. Combining this with Corollary \ref{cor:group} and using the universal coefficient theorem, one can compute the integral cohomology groups of real toric spaces $M(K,\Lambda)$, including real topological toric manifolds, in terms of $K$ and $\Lambda$. 
\begin{prop}\label{prop:group}
With the same hypotheses as Theorem \ref{thm:group}, the integral cohomology group of $M(K,\Lambda)$ is completely determined by 
\[\GG^*(K,\Lambda)=\bigoplus_{\omega\in \row\Lambda}\w H^{*-1}(K_\omega),\ \ *\geq 0,\] and the $h$-vector of $K$.	
\end{prop}
As an application of the preceding results, we get a criterion for a real topological toric manifold to be a $\mathrm{spin}^c$ manifold.
\begin{prop}\label{prop:spinc}
Let $M(K,\Lambda)$ be a real topological toric manifold. Then $M(K,\Lambda)$ is a  $\mathrm{spin}^c$ manifold if and only if $[m]\in\row\Lambda$ and the element \[w_2=\sum_{1\leq i<j\leq m}x_ix_j\in H^2(M(K,\Lambda);\Zz/2)\]
 belongs to the subgroup 
\[\Zz/2\{P_\alpha\mid\omega\in\row\Lambda,\,\alpha\in Z^1(K_\omega)\}\subset H^2(M(K,\Lambda);\Zz/2),\]
where $Z^i(K_\omega)$ is the subgroup of $\w C^i(K_\omega)$ consisting of $i$-cocycles, and $P_\alpha$ is the element defined in Theorem \ref{thm:mod 2 image}. 
\end{prop}
\begin{rem}
	Here we should mention that if the integral cohomology of a real topological toric manifold $M(K,\Lambda)$ has no elements of order $2^k$ for $k>1$, i.e., the $2$-primary component of $\GG^*(K,\Lambda)$ is zero, then whether $M(K,\Lambda)$ admits a $\mathrm{spin}^c$ structure can be read off directly from its mod $2$ cohomology ring by the rule
	\[w_3=\sum_{1\leq i<j<k\leq m}x_ix_jx_k=0\Longleftrightarrow M(K,\Lambda)\text{ is\ } \mathrm{spin}^c.\]
	This is because in this case the mod $2$ Bockstein spectral sequence of $M(K,\Lambda)$ satisfies $E_2=E_\infty$ with the first differential $d_1=Sq^1$. Note that if $M(K,\Lambda)$ is orientable then 
	\[Sq^1(w_2)=w_1w_2+w_3=w_3,\] 
	where the first equality comes from the Wu formula. However, if $\GG^*(K,\Lambda)$ has $2$-torsion, then $E_2\neq E_\infty$. So we need the information of the higher differentials $d_r$ ($r\geq 2$) of the Bockstein spectral sequence $E_*$ to know whether $w_2$ survives to $E_\infty$. But the higher differentials in $E_*$ cannot be read off directly from the mod $2$ cohomology of $M(K,\Lambda)$.
\end{rem}

\begin{proof}[Proof of Proposition \ref{prop:spinc}]
It is known that a smooth manifold $M$ admits a $\mathrm{spin}^c$ structure if and only if $M$ is orientable and its second Stiefel-Whitney class $w_2(M)$ is the mod $2$ reduction of an integral cohomology class. By \cite[Corollary 6.8 (i)]{DJ91}, the total Stiefel-Whitney class of $M(K,\Lambda)$ is given by
\[1+w_1+w_2+\cdots=\prod_{i=1}^m(1+x_i)\in H^*(M(K,\Lambda);\Zz/2).\]
(Note that the statement of \cite[Corollary 6.8 (i)]{DJ91} is for small covers, but it can be extended to real topological toric manifolds by using the real analog of \cite[Theorem 6.5]{IFM13}.)
Hence, $M(K,\Lambda)$ is orientable if and only if $w_1=x_1+\cdots+x_m=0$, i.e., $[m]\in\row\Lambda$ by Theorem \ref{thm:DJ}. Since the map 
\[\Phi_\Lambda:H^*(\A(K,\Lambda))\to H^*(\RR(K,\Lambda))\cong H^*(M(K,\Lambda))\]
is surjective by Theorem \ref{thm:group}, the image of the mod $2$ reduction $H^*(M(K,\Lambda))\to H^*(M(K,\Lambda);\Zz/2)$ is the subgroup
\[\Zz/2\{P_\alpha\mid\omega\in\row\Lambda,\,\alpha\in Z^{*-1}(K_\omega)\}\subset H^*(M(K,\Lambda);\Zz/2)\]
by Theorem \ref{thm:mod 2 image}. From here, the proposition follows immediately.
\end{proof}

In fact, there is another way to describe the cohomology group of $M(K,\Lambda)$. Set 
\begin{gather*}
	S=H^*(M(K,\Lambda);\Zz/2)\cong \Zz/2[K]/\Theta_\Lambda,\\
	I=\im \left(H^*(M(K,\Lambda))\xr{\bmod 2} S\right).
\end{gather*}
Note that $I$ can be computed using Theorem \ref{thm:mod 2 image}.
Let $\beta_{\Zz}$ be the Bockstein homomorphism $H^{*-1}(-;\Zz/2)\to H^*(-;\Zz)$ associated to the short exact sequence
\[0\to\Zz\xr{\times 2}\Zz\to\Zz/2\to 0.\]
Then there is a canonical group isomorphism
\[\bar\beta_{\Zz}:(S/I)^{*-1}\to \MM_2^*,\ \ \MM_2=\{x\in H^*(M(K,\Lambda))\mid 2x=0\}.\]
For a simplicial complex $K$, let $B^*(K)$ be the coboundary subgroup of $\w C^*(K)$. Since the differential of $\A(K,\Lambda)$ is $2\delta$, we have
\[\A_2^i=\{x\in H^i(\A(K,\Lambda))\mid 2x=0\}\cong \bigoplus_{\omega\in\row\Lambda}B^{i-1}(K_\omega)\otimes\Zz/2.\]

We claim that there is a well-defined homomorphism 
\begin{equation}\label{eq:kappa}
	\kappa:\A_2^*\to (S/I)^{*-1},\ \ [\delta\alpha]\mapsto P_\alpha,
\end{equation}
where $\alpha=\sum_{\omega\in\row\Lambda}\alpha_\omega$, $\alpha_\omega\in\w C^*(K_\omega)$, and $P_\alpha=\sum_{\omega\in\row\Lambda}P_{\alpha_\omega}$.
To see this, note that if $[\delta\alpha]=[\delta\alpha']$ in $H^*(\A(K,\Lambda))$, then $\delta\alpha-\delta\alpha'=2\delta b$ for some $b\in\A(K,\Lambda)$. So \[z:=\alpha-\alpha'-2b\in\bigoplus_{\omega\in\row\Lambda}Z^*(K_\omega).\]
It follows that $P_z=P_\alpha-P_{\alpha'}\in I$ by Theorem \ref{thm:mod 2 image}, hence $P_\alpha=P_{\alpha'}$ in $S/I$. 

\begin{prop}\label{prop:group isomorphism}
	With the same hypotheses as Theorem \ref{thm:group}, let $J\subset\A_2$ be the kernel of the map $\kappa$ in \eqref{eq:kappa}. Then there is a graded group isomorphism induced by $\Phi_\Lambda$ (or $\Psi_\Lambda$):
	\[H^*(\A(K,\Lambda))/J\cong H^*(M(K,\Lambda)).\]
\end{prop}
\begin{proof}
	Consider the commutative diagram of short exact sequences:
	\[\begin{tikzcd}
		0\arrow[r]
		&\A_2\arrow[r]\arrow[d]
		&H^*(\A(K,\Lambda))\arrow[r]\arrow[d]
		&H^*(\A(K,\Lambda))/\A_2\arrow[r]\arrow[d]
		&0\\
		0\arrow[r]
		&\MM_2\arrow[r]
		&H^*(M(K,\Lambda))\arrow[r]
		&H^*(M(K,\Lambda))/\MM_2\arrow[r]
		&0
	\end{tikzcd}\]
	where each vertical map is the composition of $\Phi_\Lambda$ (or $\Psi_\Lambda$) with the isomorphism in Theorem \ref{thm:isomorphism} (or Theorem \ref{thm:DGA}). By Theorem \ref{thm:group}, the third vertical map is an isomorphism. So the kernels of the first and the second vertical maps are isomorphic by the Snake Lemma. 
	
	Theorem \ref{thm:group} implies that the first vertical map $\A_2\to\MM_2$ is onto. We claim that this map is the composition
	\[\A_2\xr{\kappa} (S/I)\xr[\cong]{\bar\beta_{\Zz}}\MM_2.\]
	Then the isomorphism in the proposition follows immediately.
	To see this, let 
	\begin{gather*}
		T=H^*(\A(K,\Lambda)\otimes\Zz/2)=\A(K,\Lambda)\otimes\Zz/2,\\
		 N=\im\left(H^*(\A(K,\Lambda))\xr{\bmod 2} T\right).
	\end{gather*}
	We have a commutative diagram 
	\[\xymatrix{
		(T/N)\ar[r]^-{\bar\beta_{\Zz}}_-\cong\ar[d]&\A_2\ar[d]\\
		(S/I)\ar[r]^-{\bar\beta_{\Zz}}_-\cong&\MM_2}
	\]
	Hence it suffices to show that the composition $\A_2\xr{\bar\beta_{\Zz}^{-1}}(T/N)\to(S/I)$ is $\kappa$.
	Since the differential of $\A(K,\Lambda)$ is $2\delta$, $\bar\beta_{\Zz}^{-1}([\delta\alpha])=\bar\alpha+N\in T/N$ for any $\alpha\in \A(K,\Lambda)$, where $\bar\alpha$ is the mod $2$ reduction of $\alpha$. Moreover, the map $T\to S$ sends $\bar\alpha$ to $P_\alpha$, as we will see in Section \ref{sec:equivariant}. Thus, the claim holds.
\end{proof}

Next, we turn  to the cup product structure of $H^*(M(K,\Lambda))$. First we define a product structure on the graded group $\GG^*(K,\Lambda)=\bigoplus_{\omega\in \row\Lambda}\w H^{*-1}(K_\omega)$. 

For $\omega_1,\omega_2\in\row\Lambda$, let $\omega_1'=\omega_1\setminus\omega_2$ and $\omega_2'=\omega_2\setminus\omega_1$, and let $\Delta_{\omega_1+\omega_2}^{\omega_1,\omega_2}$ be the composition
\[\w H^*(K_{\omega_1})\otimes \w H^*(K_{\omega_2})\to \w H^*(K_{\omega_1'})\otimes \w H^*(K_{\omega_2'})\to \w H^*(K_{\omega_1+\omega_2}),\]
where the first map is the natural restriction and the second map is induced by the simplicial inclusion $K_{\omega_1+\omega_2}\to K_{\omega_1'}*K_{\omega_2'}$, noting that 
\[\omega_1+\omega_2=(\omega_1\cup \omega_2)\setminus(\omega_1\cap \omega_2)=\omega_1'\sqcup\omega_2'.\]
Given cohomology classes $\alpha_1\in \w H^{i-1}(K_{\omega_1})$ and $\alpha_2\in\w H^{j-1}(K_{\omega_2})$, define their $*$-product as
\[\alpha_1*\alpha_2=\Delta_{\omega_1+\omega_2}^{\omega_1,\omega_2}(\alpha_1\otimes\alpha_2)\in \w H^{i+j-1}(K_{\omega_1+\omega_2}).\]
The $*$-product  makes $\GG^*(K,\Lambda)$ into a graded ring. Our next result gives a description of the quotient ring $H^*(M(K,\Lambda))/\MM_2$ in terms of $K$ and $\Lambda$.
\begin{thm}\label{thm:ring}
Let $K$ and $\Lambda$ be as in Theorem \ref{thm:group}.	Then there is a graded ring isomorphism:
	\[H^*(M(K,\Lambda))/\MM_2\cong \GG^*(K,\Lambda),\]
where $\GG^*(K,\Lambda)$ is equipped with the $*$-product.
\end{thm}
We can further describe the cohomology ring $H^*(M(K,\Lambda))$ by using the group isomorphism in Proposition \ref{prop:group isomorphism} and endowing $H^*(\A(K,\Lambda))/J$ with a product structure. 

Given $\omega_1,\omega_2\in\row\Lambda$ and two cocycles
\[c_i=\sum_{\sigma\in K_{\omega_i}}l_{i,\sigma}\sigma^*\in Z^{p_i-1}(K_{\omega_i}),\ \ i=1,2,\ l_{i,\sigma}\in\Zz,\]
let $\omega_1'=\omega_1\setminus\omega_2$, $\omega_2'=\omega_2\setminus\omega_1$, and let $c_i'$ be the restriction of $c_i$ to $K_{\omega_i'}$, i.e.,
\[c_i'=\sum_{\sigma\in K_{\omega_i'}}l_{i,\sigma}\sigma^*\in Z^{p_i-1}(K_{\omega_i'}),\]
and assume that $\omega_1\cap\omega_2=\{k_1<\cdots<k_s\}$.
For $1\leq j\leq s$, let
\[V_{i,j}=\{\sigma\in K_{\omega_i}\mid |\sigma|=p_i,\,k_j\in\sigma,\,\sigma\cap\{k_1,\dots,k_{j-1}\}=\emptyset\},\]
and define
\begin{gather*}
	X_{i,j}=\sum_{\sigma\in V_{i,j}}(l_{i,\sigma}\bmod 2)x_{\sigma\setminus\{k_j\}}\in S^{p_i-1},\\
	E(c_1,c_2)=\sum_{j=1}^sx_{k_j}X_{1,j}X_{2,j}\in S^{p_1+p_2-1}.
\end{gather*}
Define the $\star$-product on $H^*(\A(K,\Lambda))/J$ by the rule:
\[([c_1]+J)\star([c_2]+J)=([c_1'*c_2']+J)+\bar\kappa^{-1}(E(c_1,c_2)+I),\]
where $\bar\kappa$ is the isomorphism $\A_2/J\xr{\cong}S/I$ induced by $\kappa:\A_2\to S/I$.
\begin{thm}\label{thm:ring isomorphism}
	The isomorphism in Proposition \ref{prop:group isomorphism} becomes a graded ring isomorphism if $H^*(\A(K,\Lambda))/J$ is equipped with the $\star$-product.
\end{thm}

Finally, we propose the following conjecture for real toric spaces over  Cohen-Macaulay complexes over $\Zz/2$.
\begin{conj}
	The conclusions of Theorems \ref{thm:group}, \ref{thm:ring} and \ref{thm:ring isomorphism} also hold when $K$ is a Cohen-Macaulay complex over $\Zz/2$.
\end{conj}
The remainder of the paper is devoted to proving these results. In Section \ref{sec:m-a}, we prove Theorem \ref{thm:isomorphism} by using a DGA model for $\Rr\ZZ_K$ constructed by Cai \cite{C17}. Theorem \ref{thm:DGA} will be proved by using a different DGA model for $\Rr\ZZ_K$ constructed in Section \ref{sec:DGA}. 
In Section \ref{sec:equivariant}, we construct a double complex for each subgroup $G\subset(\Zz/2)^m$, whose total cohomology equals the $G$-equivariant cohomology of $\Rr\ZZ_K$, and we use it to prove Theorem \ref{thm:mod 2 image}. 
Following the induction strategy in \cite{CC21}, we prove Theorem \ref{thm:group} for the case when $K$ is a shellable sphere in Section \ref{sec:shell}. In Section \ref{sec:PL}, we use a result of Adiprasito and Izmestiev \cite{AI15} and some fundamental results in algebraic combinatorics to complete the proof of Theorem \ref{thm:group}. In Section \ref{sec:ring}, we study the product structures of $\GG^*(K,\Lambda)$ and $H^*(\A(K,\Lambda))/J$, and give the proofs of Theorems \ref{thm:ring} and \ref{thm:ring isomorphism}.

\section{A DGA model for $\Rr\ZZ_K$}\label{sec:m-a}
In \cite{C17}, Cai constructed a DGA $\RR_K$ for any simplicial complex $K$ such that $H^*(\RR_K)$ is isomorphic to the integral cohomology ring of $\Rr\ZZ_K$. We first briefly recall his work. 

Regarding the interval $D^1=[-1,1]$ as the simplicial complex consisting of two 0-faces $\sigma_0=\{-1\}$, $\tau_0=\{1\}$ and one $1$-face $\sigma_1=\{-1,1\}$, the simplicial chain complex $(C_*(D^1),\partial)$ is the graded $\Zz$-module $\Zz\{\sigma_0,\tau_0,\sigma_1\}$ such that $\partial \sigma_1=\tau_0-\sigma_0$ and $\partial\sigma_0=\partial\tau_0=0$.
The $m$-cube $(D^1)^m$ has a natural CW-structure with each cell an $m$-fold
product of simplices of the form above. Let $D^1_i$ be the $i$th factor of $(D^1)^m$, and denote by $\sigma_{0,i}$, $\tau_{0,i}$, $\sigma_{1,i}$ the corresponding faces of $D^1_i$. By \eqref{eq:m-a}, this CW-structure of $(D^1)^m$ induces a CW-structure on the real moment-angle complex $\Rr\ZZ_K$ such that $\Rr\ZZ_K$ is embedded in $(D^1)^m$ as a cellular subcomplex, called the \emph{(canonical) cubical decomposition} of $\Rr\ZZ_K$. 

Let $(C_*(\Rr\ZZ_K),\partial)$ be the integral cellular chain complex of $\Rr\ZZ_K$ with respect to its cubical decomposition. Then $C_*(\Rr\ZZ_K)$ is a subcomplex of $\bigotimes_{i=1}^mC_*(D_i^1)$ with the boundary operator
\[\partial(e_1\otimes\cdots\otimes e_m)=\sum_{i=1}^m(-1)^{\sum_{j=1}^{i-1}\deg e_j}e_1\otimes\cdots e_{i-1}\otimes\partial e_i\otimes e_{i+1}\otimes \cdots e_m.\]
Dually, let $C^*(\Rr\ZZ_K)=\Hom(C_*(\Rr\ZZ_K),\Zz)$. Then $C^*(\Rr\ZZ_K)$ is a quotient complex of $\bigotimes_{i=1}^mC^*(D_i^1)$. 

Note that $\bigotimes_{i=1}^mC^*(D_i^1)$ is a graded ring with the cup product $\smile$ induced by the simplicial cup product on each factor $C^*(D_i^1)$. That is,
\[(e_1^*\otimes\cdots\otimes e_m^*)\smile(f_1^*\otimes\cdots\otimes f_m^*)=(-1)^s\bigotimes_{i=1}^m e^*_i\smile f_i^*, \]
where $s=\sum_{i=1}^m\deg f_i^*\sum_{j>i}\deg e_j^*$.
Let us use the following basis for $C^*(D_i^1)$,
\begin{equation}\label{eq:basis}
	1_i=\sigma_{0,i}^*+\tau_{0,i}^*,\ \ t_i=\tau_{0,i}^*,\ \ u_i=\sigma_{1,i}^*.
\end{equation}
From the definition of cup product, one checks that 
\[
	\begin{array}{lll}
	u_i\smile u_i=t_i\smile u_i=0, &u_i\smile t_i=u_i, &t_i\smile t_i=t_i,\\
	u_i\smile 1_i=1_i\smile u_i=u_i, &t_i\smile 1_i=1_i\smile t_i=t_i, & 1_i\smile1_i=1_i.
	\end{array}
\]
For a basis cochain $\bm e=e_1^*\otimes\cdots\otimes e_m^*\in \bigotimes_{i=1}^mC^*(D_i^1)$ with $e_i^*\in\{1_i,t_i,u_i\}$, let 
\[\sigma_{\bm e}=\{i\mid e_i^*=u_i\}\ \text{ and }\ \tau_{\bm e}=\{i\mid e_i^*=t_i\}.\] 
Then one can see from the above cup product structure that there is a DGA isomorphism 
\[
	\bigotimes_{i=1}^mC^*(D_i^1)\cong \RR,\ \ \bm e\mapsto u_{\sigma_{\bm e}}t_{\tau_{\bm e}},
\]
where $\RR$ and $u_{\sigma_{\bm e}}t_{\tau_{\bm e}}$ are defined in Section \ref{sec:result}. Using this description of the ring structure of $\bigotimes_{i=1}^mC^*(D_i^1)$, Cai proved that 
\begin{prop}[{\cite[Theorem 3.1 and Proposition 3.3]{C17}}]\label{prop:Cai}
Let $\RR_K$ be the DGA defined in Section \ref{sec:result}. Then there is a DGA isomorphism 
\[C^*(\Rr\ZZ_K)\cong \RR_K,\]
which induces a graded ring isomorphism:
\[H^*(\Rr\ZZ_K)\cong H^*(\RR_K).\]
\end{prop}
It is easy to check that there is a bijective cochain map (shifting the degree up by one) of cochain complexes
\begin{equation}\label{eq:sum decomposition}
	\bigoplus_{\omega\subset[m]}\w C^*(K_\omega)\xr{\cong} \RR_K,
\end{equation}
sending $\sigma^*\in \w C^*(K_\omega)$ to $u_\sigma t_{\omega\setminus\sigma}$ (see \cite[Proposition 3.3]{C17}).

Next we recall some general facts about finite group actions. For the rest of this section, $G$ will be a finite group. Let $X$ be a $G$-CW complex. (We refer to \cite[Sec. 1.1]{AP93} for the definition of $G$-CW complexes.) The orbit space $X/G$ inherits a CW structure naturally from the $G$-CW structure of $X$, so that the quotient map $\pi:X\to X/G$ is a cellular map. Then we have the induced cellular (co)chain maps
\[\pi_*:C_*(X)\to C_*(X/G),\ \ \pi^*:C^*(X/G)\to C^*(X).\]
By definition, $\pi_*$ is surjective with kernel $\{e-g\cdot e\mid e\in C_*(X),\,g\in G\}$, so  
\[C_*(X/G)\cong \Zz\otimes_{\Zz G}C_*(X).\] 
Dually, $\pi^*$ is injective with image $C^*(X)^G$, the invariant elements of $C^*(X)$ under the $G$-action, where the action of $G$ on $C^*(X)$ is given by $(g\cdot\varphi)(e)=\varphi(g^{-1}\cdot e)$ for any $g\in G$, $\varphi\in C^*(X)$, $e\in C_*(X)$. Indeed, we have the following sequence of isomorphisms
\begin{equation}\label{eq:equivariant cohomology}
	\begin{split}
C^*(X)^G\cong \Hom_{\Zz G}(\Zz,C^*(X))&\cong\Hom(\Zz\otimes_{\Zz G}C_*(X),\Zz)\\
&\cong \Hom(C_*(X/G),\Zz)\cong C^*(X/G),
	\end{split}
\end{equation}
using \cite[Proposition II.5.2]{CE56} for the second isomorphism.

In particular, for a free $G$-CW complex $X$ (i.e., $\pi:X\to X/G$ is a cellular covering with deck transformation group $G$), we have the following standard fact.
\begin{lem}\label{lem:norm}
	Let $X$ be a free $G$-CW complex. Then $C^*(X)^G$ is the subcochain complex
	\[\{\sum_{g\in G}g\cdot \varphi\mid \varphi\in C^*(X)\}.\]
\end{lem}
\begin{proof}
	Since $X$ is a free $G$-CW complex, $C_*(X)$ is a free $\Zz G$-module, hence $C^*(X)$ is a weakly projective $\Zz G$-module by \cite[Propositions X.8.5 and XII.1.1]{CE56}. So we obtain the desired identification by \cite[Proposition XII.1.3]{CE56}.
\end{proof}
Now we turn our attention to the special case $X=\Rr\ZZ_K$ with the cubical decomposition and $G=\ker\Lambda$ for some characteristic function $\Lambda$ over $K$. One checks that for every $g\in\ker\Lambda$ and every cell $e\subset\Rr\ZZ_K$, $g\cdot e$ is again a cell of $\Rr\ZZ_K$.  Moreover, since $\ker\Lambda$ acts freely on $\Rr\ZZ_K$, $g\cdot e\neq e$ for any cell $e\in \Rr\ZZ_K$ and any $g\neq 1$ by the Brouwer fixed point theorem. Hence $\Rr\ZZ_K$ is a free $G$-CW-complex.

Let us consider the action of $\ker\Lambda$ on $C^*(\Rr\ZZ_K)$. Let $r:D^1\to D^1$, $x\mapsto -x$, be the reflection. Then the induced map $r^*$ on $C^*(D^1)$ is given by
\[r^*(\sigma_0^*)=\tau_0^*,\ \ r^*(\tau_0^*)=\sigma_0^*,\ \ r^*(\sigma_1^*)=-\sigma_1^*.\]
It follows that for $\bm a_i$ the $i$th coordinate vector of $(\Zz/2)^m$, the action of $\bm a_i$ on the basis \eqref{eq:basis} is given by
\[\bm a_i\cdot 1_j=1_j,\ \ 
\bm a_i\cdot t_j=\begin{cases}
	t_j, &j\neq i,\\
	1_j-t_j, &j=i,
\end{cases}\ \  
\bm a_i\cdot u_j=\begin{cases}
u_j, &j\neq i,\\
-u_j, &j=i.
\end{cases}\]
Therefore, under the isomorphism \eqref{eq:sum decomposition}, one easily checks that the action of an element $g\in\ker\Lambda$ on $\RR_K$ is given by the map $\phi_g$ defined in Section \ref{sec:result}. 

We can now prove Theorem \ref{thm:isomorphism}.
\begin{proof}[Proof of Theorem \ref{thm:isomorphism}]
	By Proposition \ref{prop:Cai}, $C^*(\Rr\ZZ_K)\cong \RR_K$. Then \eqref{eq:equivariant cohomology} shows that
	 \[C^*(M(K,\Lambda))\cong (\RR_K)^{\ker\Lambda}=\RR(K,\Lambda),\]
	 and so $C^*(M(K,\Lambda))\otimes\kk\cong \RR(K,\Lambda)\otimes\kk$, which gives the second assertion of the theorem. The first assertion follows from Lemma \ref{lem:norm}. 
\end{proof}
We end this section by proving the following result mentioned in Section \ref{sec:result}.
\begin{lem}\label{lem:invariants}
	Let $\Lambda:(\Zz/2)^m\to (\Zz/2)^n$ be a surjective linear map. Then for any $\omega\in\row\Lambda$ and any oriented face $\sigma\in K_\omega$, we have
	$f_{\sigma,\omega}\in \RR(K,\Lambda)$ and $h_{\sigma,\omega}\in \QQ(K,\Lambda)$, where $f_{\sigma,\omega},h_{\sigma,\omega}$ are the elements defined in \eqref{eq:f,h}. Moreover, the maps $\varphi_\omega$ and $\psi_\omega$ defined in \eqref{eq:map phi_omega} are chain maps.
\end{lem}
\begin{proof}
	We prove only the results for $f_{\sigma,\omega}$, since the argument for $h_{\sigma,\omega}$ is very similar. By \eqref{eq:action}, $g\cdot f_{\sigma,\omega}$ ($g\in\ker\Lambda$) is a linear combination 
	\begin{equation}\label{eq:linear}
		\sum_{\tau\subset\omega\setminus\sigma}n_\tau u_\sigma t_\tau,\ \ n_\tau\in\Zz.
	\end{equation}
	We need to show that $n_\tau=(-2)^{|\tau|}$ for any $g\in\ker\Lambda$. 
		
	For an element $g\in\ker\Lambda$ and a monomial $u_\sigma t_{\omega'}$ with $\omega'\subset\omega\setminus\sigma$, $g\cdot u_\sigma t_{\omega'}$ contributes a nonzero term $\pm u_\sigma t_\tau$ to \eqref{eq:linear} only if $\tau\subset\omega'\subset g\cup\tau$, by \eqref{eq:action} again. Hence, $n_\tau$ is the coefficient of $u_\sigma t_\tau$ in the linear combination
	\begin{equation}\label{eq:linear2}
	g\cdot\sum_{\tau\subset\omega'\subset (g\cup\tau)\cap (\omega\setminus\sigma)}(-2)^{|\omega'|}u_\sigma t_{\omega'}.
    \end{equation}
	For each $\tau\subset\omega\setminus\sigma$, let $k_\tau=|g\cap(\omega\setminus(\sigma\cup \tau))|$, $l_\tau=|g\cap(\sigma\cup\tau)|$. Then there are exactly $\binom{k_\tau}{i}$ monomials $u_\sigma t_{\omega'}$ such that $\tau\subset\omega'\subset(g\cup\tau)\cap (\omega\setminus\sigma)$ and $|\omega'\setminus\tau|=i$, each of which contributes a monomial term $(-1)^{l_\tau}u_\sigma t_\tau$ to \eqref{eq:linear2}. Adding them together, one computes that the coefficient of $u_\sigma t_\tau$ in \eqref{eq:linear2} is
	\[(-1)^{l_\tau}\sum_{i=0}^{k_\tau}\binom{k_\tau}{i}(-2)^{|\tau|+i}=(-1)^{l_\tau}(-2)^{|\tau|}(1-2)^{k_\tau}=(-2)^{|\tau|}(-1)^{k_\tau+l_\tau}.\]
	Since $k_\tau+l_\tau=|g\cap\omega|$ is even by \cite[Lemma 3.2]{CP20}, we get $n_\tau=(-2)^{|\tau|}$.
	
	An easy calculation shows that $2\varphi_\omega\delta(\sigma^*)=d\varphi_\omega(\sigma^*)$ for any $\sigma\in K_\omega$, proving the final assertion of the lemma.
\end{proof}

\section{A DGA model for $M(K,\Lambda)$}\label{sec:DGA}
In this section, we assume that the coefficient ring is $\Zz$. All the results can be readily extended to the case of an arbitrary coefficient ring $\kk$. 

We will use the notation $S_*(X)$ and $S^*(X)$ for the normalized singular chain and cochain complexes of a space $X$ respectively (see \cite[Sec. I.4.(a)]{FHT01} for the definition). $S_*(X)$ is a DGC (differential graded coalgebra) with the Alexander-Whitney diagonal $\Delta$ as the coproduct:
\[\Delta_X:S_*(X)\to S_*(X)\otimes S_*(X),
\ \ \Delta_X(\sigma)=\sum_{i=0}^n\sigma|_{[v_0,\ldots,v_i]}\otimes\sigma|_{[v_i,\ldots,v_n]}, 
\]
where $\sigma:\Delta^n\to X$ is a singular $n$-simplex (the vertex set of $\Delta^n$ is $\{v_0,\ldots,v_n\}$) and $\sigma|_{[v_i,\dots,v_j]}$ ($i\leq j$) is the restriction of $\sigma$ to the face spanned by $\{v_i,v_{i+1},\ldots,v_j\}$.
Dually, $S^*(X)$ is a DGA with the usual cup product. If a finite group $G$ acts on $X$, then $S_*(X)$ becomes a $\Zz G$-DGC and $S^*(X)$ a $\Zz G$-DGA.

As we mentioned in Section \ref{sec:result}, $\RR(K,\Lambda)$ is not closed under multiplication. The reason for this is that the inclusion $C_*(D^1)\subset S_*(D^1)$ is not $\Zz/2$-equivariant. So we start by constructing a larger subchain complex $T_*(D^1)$ of $S_*(D^1)$ so that the inclusion $T_*(D^1)\subset S_*(D^1)$ is a $\Zz/2$-equivariant quasi-isomorphism.

Let $\Delta^n$ be the $n$-simplex with vertices $v_0,\dots,v_n$, and let $\sigma_n:\Delta^n\to D^1$ (resp. $\tau_n:\Delta^n\to D^1$) be the simplicial map  sending $v_i$ to $(-1)^{i+1}$ (resp. to $(-1)^i$). Define $T_*(D^1)$ as the submodule of $S_*(D^1)$ generated by $\sigma_n,\tau_n$, $n\geq 0$.  $T_*(D^1)$ is clearly a chain complex, and the boundary maps are given by
\[\partial \sigma_n=\tau_{n-1}+(-1)^n\sigma_{n-1},\quad \partial \tau_n=\sigma_{n-1}+(-1)^n\tau_{n-1}.\]

\begin{lem}\label{lem:quasi-iso}
	The inclusion $i:T_*(D^1)\to S_*(D^1)$ is a $\Zz/2$-equivariant DGC map and a quasi-isomorphism.
\end{lem}
\begin{proof}
	It is straightforward to check that $i$ is a DGC map and a quasi-isomorphism. 
	Let $r:D^1\to D^1$ be the reflection as before. Then by definition, we have
	\[r\sigma_n=\tau_n,\quad r\tau_n=\sigma_n,\]
	hence $i$ is $\Zz/2$-equivariant.
\end{proof}

For each face $\sigma\in K$, let 
\[Y_\sigma=\prod_{i=1}^m Y_{\sigma|i},\ \text{ where }\ Y_{\sigma|i}=\begin{cases}D^1&\text{if } i\in\sigma,\\
	S^0&\text{if }i\not\in\sigma.\end{cases}\]
We define
\[S(\sigma)=\bigotimes_{i=1}^m S_*(Y_{\sigma|i})\ \text{ and }\ S(K)=\colim_{\sigma\in K}S(\sigma),\]
where the colimit is defined by the inclusions $S(\sigma)\subset S(\tau)$ with $\sigma\subset\tau\in K$. Then $S(K)$ is a subchain complex of $\bigotimes_{i=1}^m S_*(D^1_i)$. Dually, define
\[S(K)^*=\Hom(S(K),\Zz)\] 
is a quotient DGA of $\bigotimes_{i=1}^m S^*(D^1_i)$. Similarly, set $T_*(S^0)=S_*(S^0)$ and define
\begin{gather*}
	T(\sigma)=\bigotimes_{i=1}^m T_*(Y_{\sigma|i}),\ \ T(K)=\colim_{\sigma\in K}T(\sigma),\\
	T(K)^*=\Hom(T(K),\Zz).
\end{gather*}

\begin{prop}\label{prop:DGC model}
	The inclusion $T(K)\subset S(K)$ is a $(\Zz/2)^m$-equivariant DGC map and a quasi-isomorphism.
\end{prop}
\begin{proof}
The map is clearly a $(\Zz/2)^m$-equivariant DGC map, so we only verify that it is a quasi-isomorphism. If $K$ is a simplex, then the result follows by Lemma \ref{lem:quasi-iso}. For the general case, we use induction on the number of faces in $K$.	Suppose that $K$ is obtained from $K'$ by attaching a simplex $2^\sigma$. Then there is a commutative diagram with exact rows
	\[\xymatrix{
	0\ar[r]&T(K'\cap \sigma)\ar[r]\ar[d]&T(K')\oplus T(\sigma)\ar[r]\ar[d]&T(K)\ar[r]\ar[d]&0\\
	0\ar[r]&S(K'\cap \sigma)\ar[r]&S(K')\oplus S(\sigma)\ar[r]&S(K)\ar[r]&0}
\]
By induction, the first two vertical maps are quasi-isomorphisms. Then, applying the five lemma to the induced long exact sequences of cohomology, we get the desired isomorphism $H^*(T(K))\cong H^*(S(K))$.
\end{proof}
This proof basically follows the argument used in \cite[Proposition 3.5]{Franz22}. The following proposition is the dual version of Proposition \ref{prop:DGC model}.
\begin{prop}\label{prop:DGA model}
	The map $S(K)^*\to T(K)^*$ is a $(\Zz/2)^m$-equivariant DGA map and a quasi-isomorphism.
\end{prop}

Recall that for two spaces $X$ and $Y$, there is the Eilenberg-Zilber map
\[\nabla: S_*(X)\otimes S_*(Y)\to S_*(X\times Y),\]
which is a chain homotopy equivalence, natural in the pair $(X,Y)$, and is a DGC map (\cite[(17.6)]{EM66}). Let 
\[\nabla^*:S^*(X\times Y)\to \Hom(S_*(X)\otimes S_*(Y),\Zz)\] 
be the dual of $\nabla$.  By iteration, the Eilenberg-Zilber map and its dual extend to more than two factors. 

For each $\sigma\in K$, we have the Eilenberg-Zilber map $\nabla_\sigma:S(\sigma)\to S_*(Y_\sigma)$, and all $\nabla_\sigma$ fit together to give a chain map $\nabla_K:S(K)\to S_*(\Rr\ZZ_K)$. Dually, we have the cochain map $\nabla_K^*:S^*(\Rr\ZZ_K)\to S(K)^*$.
The following two propositions are analogs of Proposition \ref{prop:DGC model} and \ref{prop:DGA model}.
\begin{prop}
The map $\nabla_K:S(K)\to S_*(\Rr\ZZ_K)$ is a $(\Zz/2)^m$-equivariant DGC map and  a chain homotopy equivalence.
\end{prop}
\begin{proof}
The proof is similar to the proof of Proposition \ref{prop:DGC model}. For the induction step, suppose that $K=K'\cup 2^\sigma$.  Consider the commutative diagram with exact rows
\[\xymatrix{
	0\ar[r]&S(K'\cap\sigma)\ar[r]\ar@{<->}[d]^{\simeq}&S(K')\oplus S(\sigma)\ar[r]\ar@{<->}[d]^{\simeq}&S(K)\ar[r]\ar@{<->}[d]&0\\
	0\ar[r]&S_*(\ZZ_{K'\cap\sigma})\ar[r]&S_*(\Rr\ZZ_{K'})\oplus S_*(\Rr\ZZ_\sigma)\ar[r]&S_*(\Rr\ZZ_{K'})+S_*(\Rr\ZZ_\sigma)\ar[r]&0}
\]
By induction, the first two vertical maps are chain homotopy equivalences.
Hence the third vertical map is also a chain homotopy equivalence.
Observe that $\nabla_K$ factors as 
\[S(K)\to S_*(\Rr\ZZ_{K'})+S_*(\Rr\ZZ_\sigma)\to S_*(\Rr\ZZ_K),\] 
where the second map is a chain homotopy equivalence by excision.
\end{proof}
\begin{prop}\label{cor:E-Z}
	The map $\nabla_K^*:S^*(\Rr\ZZ_K)\to S(K)^*$ is a $(\Zz/2)^m$-equivariant DGA map and a cochain homotopy equivalence.
\end{prop}

In fact, $T^*(D^1)$, $\bigotimes_{i=1}^m T^*(D^1_i)$, $T(K)^*$ are just the DGAs $\Gamma$, $\QQ$, $\QQ_K$ defined in the discussion preceding Theorem \ref{thm:DGA}, respectively. To see this, let $\sigma_n^*$, $\tau_n^*\in T^*(D^1)$ be the duals of $\sigma_n$, $\tau_n$ respectively.
Then there is an obvious $\Zz/2$-equivariant DGA isomorphism $T^*(D^1)\cong\Gamma$ given by $\sigma_n^*\mapsto \mu_n$, $\tau_n^*\mapsto \mu_n'$, which induces a $(\Zz/2)^m$-equivariant DGA isomorphism  $\bigotimes_{i=1}^m T^*(D^1_i)\cong\QQ$, and it passes to the quotient to give a DGA isomorphism $T(K)^*\cong\QQ_K$.

As in the proof of Lemma \ref{lem:norm}, the first assertion of Theorem \ref{thm:DGA} follows from \cite[Propositions X.8.5, XII.1.1 and XII.1.3]{CE56}, since $T(K)$ is a free $\Zz G$-module for $G=\ker\Lambda$ when $\Lambda$ is a characteristic function over $K$. 

Recall from \eqref{eq:equivariant cohomology} that for a finite group $G$ acting on a space $X$, we have a DGA isomorphism $S^*(X/G)\cong S^*(X)^G$, so the second assertion of Theorem \ref{thm:DGA} is a consequence of the following result.
\begin{thm}\label{thm:quasi-iso}
	Let $\Lambda:(\Zz/2)^m\to (\Zz/2)^n$ be a characteristic function over $K$. Then the DGA maps 
	\[S^*(\Rr\ZZ_K)^{\ker\Lambda}\xr{\nabla_K^*}(S(K)^*)^{\ker\Lambda}\ \text{ and }\  (S(K)^*)^{\ker\Lambda}\to (T(K)^*)^{\ker\Lambda}\]
	 are  quasi-isomorphisms.
\end{thm}
Since $S_*(\Rr\ZZ_K)$, $S(K)$, $T(K)$ are graded free $\Zz G$-modules, by \cite[Proposition X.8.5]{CE56}, $S^*(\Rr\ZZ_K)$, $S(K)^*$, $T(K)^*$  are graded weakly injective $\Zz G$-modules, i.e., they are weakly injective in each degree (see \cite[Sec. X.8]{CE56} for the definition and properties of a weakly injective module). Hence, Theorem \ref{thm:quasi-iso} is a consequence of the following general result in homological algebra, using Propositions \ref{prop:DGA model} and \ref{cor:E-Z}. 
\begin{prop}\label{prop:weakly injective}
Let $G$ be a group, and assume that $A_1,A_2$ are two graded weakly injective $\Zz G$-modules satisfying:
\begin{itemize}
	\item $A_1,A_2$ are DGAs over $\Zz$ (with differentials of degree $1$, both denoted by $\delta$), and for any $g\in G$, the maps $A_1\xr{g\cdot}A_1$ and $A_2\xr{g\cdot}A_2$ are DGA maps;
	\item There is a $G$-equivariant DGA map $A_1\to A_2$ that is a quasi-isomorphism.
\end{itemize}	
Then the restriction DGA map $A_1^G\to A_2^G$ is a quasi-isomorphism.
\end{prop}
\begin{proof}
Let $\PP:\cdots\xr{\partial} P_k\xr{\partial}\cdots \xr{\partial}P_0\to\Zz$ be a projective resolution of the trivial $\Zz G$-module $\Zz$. Then for $i=1,2$, the graded $\Zz$-module
\begin{gather*}
\DD_i=\Hom_{\Zz G}(\PP,A_i)=\bigoplus_{j,k}\Hom_{\Zz G}(P_k,A_i^j)\\
\DD_i^n=\bigoplus_{j+k=n}\Hom_{\Zz G}(P_k,A_i^j)
\end{gather*}
 is a double complex with total differentials  given by
\begin{equation}\label{eq:diff}
	d\varphi=\delta\varphi+(-1)^j\varphi\partial\ \text{ for }\ \varphi\in \Hom_{\Zz G}(P_k,A_i^j).
\end{equation}

Filtering $\DD_i$, $i=1,2$, by the degree of $A_i$, we get a spectral sequence $^\mathrm{I}E_*^{*,*}(\DD_i)$, whose $E_1$-term is
\[^\mathrm{I}E_1^{p,q}(\DD_i)=H^q(\Hom_{\Zz G}(\PP,A_i^p))=\Ext^q_{\Zz G}(\Zz,A_i^p).\]
Since $A_i^p$ are weakly injective by hypothesis, $\Ext^q_{\Zz G}(\Zz,A_i^p)=0$ for $q>0$ (see \cite[Proposition X.8.2a]{CE56}). So $^\mathrm{I}E_*^{*,*}(\DD_i)$ collapses: 
\begin{equation}\label{eq:E_1}
	^\mathrm{I}E_1^{*,*}(\DD_i)={^\mathrm{I}}E_1^{*,0}(\DD_i)=\mathrm{Ext}^0_{\Zz G}(\Zz,A_i)=A_i^G.
\end{equation}
Hence we have an isomorphism of $\Zz$-modules 
\begin{equation}\label{eq:E_2}
	H^*(\DD_i,d)\cong H^*(A_i^G,\delta).
\end{equation}
On the other hand, filtering $\DD_i$ by the degree of $\PP$, we get another spectral sequence $^\mathrm{II}E_*^{*,*}(\DD_i)$, whose $E_1$-term is
\begin{equation}\label{eq:E_1'}
	^\mathrm{II}E_1^{p,q}(\DD_i)=H^q(\Hom_{\Zz G}(P_p,A_i))=\Hom_{\Zz G}(P_p,H^q(A_i)),
\end{equation}
where the second equality comes from the fact that $P_p$ is projective, so $\Hom_{\Zz G}(P_p,-)$ is an exact functor.

Since the map $A_1\to A_2$ is $G$-equivariant, it induces a map of double complexes
$\DD_1\to\DD_2$, which preserves filtrations. Hence, there are induced spectral sequence maps
\[^\mathrm{I}E_*^{*,*}(\DD_1)\to {^\mathrm{I}E}_*^{*,*}(\DD_2)\ \text{ and }\ ^\mathrm{II}E_*^{*,*}(\DD_1)\to {^\mathrm{II}E}_*^{*,*}(\DD_2).\]
By \eqref{eq:E_1'}, the second map on the $E_1$-page is an isomorphism, using
the assumption that $A_1\to A_2$ is a quasi-isomorphism. An application of the spectral sequence comparison theorem (see, for example, \cite[Theorem 5.2.12]{Wei94}) shows that $\DD_1\to\DD_2$ is a quasi-isomorphism. On the other hand, by \eqref{eq:E_1}, the first map on the $E_1$-page is just the DGA map $A_1^G\to A_2^G$. Then the proposition follows from \eqref{eq:E_2}.
\end{proof}

\section{The equivariant cohomology of $\Rr\ZZ_K$}\label{sec:equivariant}
Given a (discrete) group $G$, there is the universal principal $G$-bundle $EG\to BG$ with the total space $EG$ contractible. Let $X$ be a $G$-space. The diagonal $G$-action on $EG\times X$, given by
\[g(e,x)=(ge,gx),\ \ g\in G,\ e\in EG,\ x\in X,\]
is free. Its orbit space is denoted by $EG\times_GX$ and is called the \emph{Borel construction} of $X$ by $G$. The ordinary cohomology $H^*(EG\times_GX;\kk)$ is called the \emph{$G$-equivariant cohomology} of $X$ and is denoted $H^*_G(X;\kk)$. In this section, we construct a double complex $\DD(G,K)$ for each subgroup $G\subset(\Zz/2)^m$ so that the total cohomology group of $\DD(G,K)\otimes\kk$ is isomorphic to $H^*_G(\Rr\ZZ_K;\kk)$.

Write the cyclic group $C_2=\Zz/2$ multiplicatively with generator $g$ and identity element $1$. Then there is a simple $\Zz C_2$-free resolution  of $\Zz$
\[F_*:\cdots\xrightarrow{\cdot(g-1)}F_2\xrightarrow{\cdot(g+1)}F_1\xrightarrow{\cdot(g-1)}F_0\to\Zz,\]  
where $F_i=\Zz C_2$, $i\geq 0$, with differentials \[F_{2j+1}\xrightarrow{\cdot(g-1)}F_{2j}\ \text{ and }\ F_{2j}\xrightarrow{\cdot(g+1)}F_{2j-1}.\] 
Thus, for $G=(\Zz/2)^n$ the $n$-fold tensor product $\PP_G=F_*^{\otimes n}$ is a $\Zz G$-free  resolution of $\Zz$.
It has the form  \[\PP_G=\bigoplus_{i\geq0}P_{G,i}=\Zz[n]\otimes\Zz G,\]
where $\Zz[n]=\Zz[x_1,\dots,x_n]$ is the graded polynomial ring with each $x_i$ of degree $1$ and $P_{G,i}=\Zz[n]_i\otimes\Zz G$. We write $g_i$ for the generator of the $i$th summand $\Zz/2$ of $G$. Then the differential of $\PP_G$ is given by
\begin{equation}\label{eq:differential of P}
	\partial(x_1^{k_1}\cdots x_n^{k_n}\otimes f)=\sum_{i=1}^n(-1)^{\sum_{j=1}^{i-1}k_j}x_1^{k_1}\cdots x_i^{k_i-1}\cdots x_n^{k_n}\otimes\varepsilon_if,
\end{equation}
where 
\[\varepsilon_i=\begin{cases}
	g_i+1,&\text{if $k_i\neq 0$ is even,}\\
	g_i-1,&\text{if $k_i$ is odd,}\\
	0,&\text{if }k_i=0.
\end{cases}\]

Given a subgroup $G\subset(\Zz/2)^m$ with $G\cong(\Zz/2)^n$ for some $n\leq m$, we have the following $\Zz$-module isomorphisms
\begin{equation}\label{eq:adjoint isom}
	\begin{split}
		&\Hom_{\Zz G}(\PP_G,C^*(\Rr\ZZ_K))=\bigoplus_{i,j}\Hom_{\Zz G}(P_{G,i},C^j(\Rr\ZZ_K))\\
		\cong&\bigoplus_{i,j}\Hom(\Zz[n]_i,\Hom_{\Zz G}(\Zz G,C^j(\Rr\Zz_K))) \\
		\cong&\bigoplus_{i,j}\Hom(\Zz[n]_i,C^j(\Rr\ZZ_K))\\
		\cong&\bigoplus_{i,j}\Hom(\Zz[n]_i\otimes C_j(\Rr\ZZ_K),\Zz)\\
		\cong&\Hom(\Zz[n],\Zz)\otimes C^*(\Rr\ZZ_K)\cong\Zz[n]\otimes\RR_K,
	\end{split}
\end{equation}
where the first and third isomorphisms come from the Hom-tensor adjunction. We use the differentials as in \eqref{eq:diff} to make $\Hom_{\Zz G}(\PP_G,C^*(\Rr\Zz_K))$ into a double complex. Then, using \eqref{eq:adjoint isom}, $\Zz[n]\otimes\RR_K$ becomes a 
double complex with total differentials given by
\begin{equation}\label{eq:differential D(G,K)}
	\begin{split}
		d(x_1^{k_1}&\cdots x_n^{k_n}\otimes\alpha)=\sum_{i=1}^n(-1)^{\sum_{j=1}^{i-1}k_j}x_1^{k_1}\cdots x_i^{k_i+1}\cdots x_n^{k_n}\otimes\varepsilon_i\alpha\\
		&+(-1)^{\sum_{j=1}^nk_j}x_1^{k_1}\cdots x_n^{k_n}\otimes d\alpha
		,\end{split}
\end{equation}
where 
\[\varepsilon_i=\begin{cases}
	g_i-1,&\text{if $k_i$ is even,}\\
	g_i+1,&\text{if $k_i$ is odd.}\\
\end{cases}\]
Here we think of $g_i$ as an element of $(\Zz/2)^m$. 

Let $\DD(G,K)=\Zz[n]\otimes\RR_K$ with the above double complex structure. Then $\DD(G,K)$ is functorial with respect to inclusions of subgroups of $(\Zz/2)^m$ and inclusions of subcomplexes (viewed as simplicial complexes on $[m]$) of the $(m-1)$-simplex $2^{[m]}$ in the following sense. Given inclusions $(\Zz/2)^n\cong G\subset G'\cong (\Zz/2)^{n'}$ and $K\subset K'$, since $G$ and $G'$ are vector spaces over $\Zz/2$ we can choose a basis $\{g_1,\ldots,g_n\}$ for $G$ and extend it to a basis $\{g_1,\ldots,g_n,g_{n+1},\ldots,g_{n'}\}$ for $G'$. Then we have induced maps
\begin{equation}\label{eq:induced maps}
	 \phi:\Zz[n']\to\Zz[n],\quad \psi:\RR_{K'}\to \RR_K,
\end{equation}
where $\phi$ is the ring homomorphism given by $\phi(x_i)=x_i$ for $1\leq i\leq n$ and $\phi(x_i)=0$ for $i>n$. One checks that the tensor product map
\[\phi\otimes\psi:\DD(G',K')\to \DD(G,K)\]
is a double complex map with the differential given in \eqref{eq:differential D(G,K)}, where $\varepsilon_i$ is defined using the above basis extension. 

Notice that this double complex map is not canonical because it depends on the choice of basis extension. In fact, there is another slightly more complicated  double complex that is canonically functorial, and we will discuss it later.

\begin{prop}\label{prop:equivariant}
	There is an isomorphism of graded $\kk$-modules for any subgroup $G\subset (\Zz/2)^m$ and any coefficient ring $\kk$: \[H^*(\DD(G,K)\otimes\kk)\cong H^*_G(\Rr\ZZ_K;\kk).\]
	The isomorphism is natural with respect to the inclusions of subgroups of $(\Zz/2)^m$ and subcomplexes of $2^{[m]}$ in the sense described above. 
\end{prop}

Before proving the proposition, let us recall the canonical cellular decomposition of $E\Zz/2=S^\infty=\colim_{n\geq 0}S^n$. Inductively, each sphere $S^n$ can be obtained from
$S^{n-1}$ by attaching two cells of dimension $n$:
\[S^n=D^n\cup_fS^{n-1}\cup_hD^n.\]
Here the attaching maps $f,h:\partial D^n\to S^{n-1}$ are the identity. The group $\Zz/2$ freely acts on $S^\infty$ by the antipodal map that is cellular, so that \[B\Zz/2=S^\infty/(\Zz/2)=\Rr P^\infty.\]

For a subgroup $G\cong (\Zz/2)^n\subset(\Zz/2)^m$, let $C^*(EG\times\Rr\ZZ_K)$ be the cellular cochain complex with respect to the canonical cellular structure on the product space $EG\times\Rr\ZZ_K=(S^\infty)^{\times n}\times \Rr\ZZ_K$. Then the $G$-action on $EG\times\Rr\ZZ_K$ is a free cellular action, and the projection $EG\times\Rr\ZZ_K\to \Rr\ZZ_K$ is a $G$-equivariant cellular map.
\begin{proof}[Proof of Proposition \ref{prop:equivariant}]
	Suppose that $G\cong(\Zz/2)^n$, and let $\PP_G(\kk)=\Zz[n]\otimes\kk G$ be the $\kk G$-free resolution of $\kk$. Since $C_*(EG\times\Rr\ZZ_K;\kk)$ is a free $\kk G$-module, the cellular cochain complex $C^*(EG\times\Rr\ZZ_K;\kk)$ is a weakly injective $\kk G$-module by \cite[Proposition X.8.5]{CE56}. Hence the cochain map
	\[\Hom_{\kk G}(\kk,C^*(EG\times\Rr\ZZ_K;\kk))\to \Hom_{\kk G}(\PP_G(\kk),C^*(EG\times\Rr\ZZ_K;\kk))\]
	is a quasi-isomorphism (see the proof of \eqref{eq:E_2}). Note that \[C^*(EG\times\Rr\ZZ_K;\kk)^G\cong\Hom_{\kk G}(\kk,C^*(EG\times\Rr\ZZ_K;\kk)),\]
	and that \[H^*(C^*(EG\times\Rr\ZZ_K;\kk)^G)\cong H^*_G(\Rr\ZZ_K;\kk),\]
	using \eqref{eq:equivariant cohomology}.
	
	If we filter the above double complex $\Hom_{\kk G}(\PP_G(\kk),C^*(EG\times\Rr\ZZ_K;\kk))$ by the degree of $\PP_G(\kk)$, then the $E_1$-term of the corresponding spectral sequence is of the form
	\[E_1^{p,q}=\Hom_{\kk G}(P_{G,p}\otimes\kk,H^q(EG\times\Rr\ZZ_K;\kk))\]
	(see the proof of \eqref{eq:E_1'}). Let $\pi:EG\times\Rr\ZZ_K\to \Rr\ZZ_K$ be the projection. Then there is an induced double complex map 
	\[\pi^*:\Hom_{\kk G}(\PP_G(\kk),C^*(\Rr\ZZ_K;\kk))\to \Hom_{\kk G}(\PP_G(\kk),C^*(EG\times\Rr\ZZ_K;\kk)),\]
	and a spectral sequence map which on the $E_1$-page is
	\[\pi^*:\Hom_{\kk G}(P_{G,p}\otimes\kk,H^q(\Rr\ZZ_K;\kk))\to \Hom_{\kk G}(P_{G,p}\otimes\kk,H^q(EG\times\Rr\ZZ_K;\kk)).\]
	Since $\pi^*: H^*(\Rr\ZZ_K;\kk)\to H^*(EG\times\Rr\ZZ_K;\kk)$ is a $G$-equivariant isomorphism, it follows that $\Hom_{\kk G}(\PP_G(\kk),C^*(EG\times\Rr\ZZ_K;\kk))$ is quasi-isomorphic to 
	\begin{equation}\label{eq:isomorphism D(K,G)}
		\Hom_{\kk G}(\PP_G(\kk),C^*(\Rr\ZZ_K;\kk))\cong\DD(G,K)\otimes\kk,
\end{equation}	
	using the spectral sequence comparison theorem and \eqref{eq:adjoint isom}.
	From here, the first assertion of the proposition follows immediately. 
	
	For naturality, suppose that we have inclusions $(\Zz/2)^n\cong G\subset G'\cong (\Zz/2)^{n'}\subset (\Zz/2)^m$ and $K\subset K'$ , where $K$ and $K'$ are simplicial complexes on $[m]$. We use the convention that there is a basis extension of $G\subset G'$ which gives the cellular inclusion 
	\[\begin{split}
		EG\times \Rr\ZZ_K=(S^\infty)^{\times n}\times \Rr\ZZ_K&\subset (S^\infty)^{\times n'}\times \Rr\ZZ_{K'}=EG'\times \Rr\ZZ_{K'}\\
		(y_1,\ldots,y_n,x)&\mapsto(y_1,\ldots,y_n,*\ldots,*,x).
	\end{split}\]
	Here $*\in S^0\subset S^\infty$ is a fixed point. Then we have the composition of cochain maps
	\[C^*(EG'\times\Rr\ZZ_{K'};\kk)^{G'}\hookrightarrow C^*(EG'\times\Rr\ZZ_{K'};\kk)^G\to C^*(EG\times\Rr\ZZ_K;\kk)^G,\] 
	which induces the homomorphism of equivariant cohomology
	\[H^*_{G'}(\Rr\ZZ_{K'};\kk)\to H^*_G(\Rr\ZZ_K;\kk).\]
	Furthermore, we see that the quasi-isomorphisms in the previous arguments are all functorial with respect to the above cellular inclusion and the maps in \eqref{eq:induced maps}. This completes the proof of the proposition.
\end{proof}

To prove Theorem \ref{thm:mod 2 image}, let us construct another double complex to compute the equivariant cohomology of $\Rr\ZZ_K$. Let $H=(\Zz/2)^m$ with the standard coordinate basis $\{g_1,\ldots,g_m\}$ and let $\PP=\PP_H=\Zz[m]\otimes\Zz H$ be the chain complex with differentials given by \eqref{eq:differential of P}. Note that for any subgroup $G\subset H$, $\Zz H$ is a $\Zz G$-free module, so $\PP$ can be viewed as a $\Zz G$-free resolution of $\Zz$. For any subgroup $G\subset H$ and any simplicial complex $K$ on $[m]$, define
\[\PP(G,K)=\Hom_{\Zz G}(\PP,C^*(\Rr\ZZ_K)).\]
$\PP(G,K)$ becomes a double complex with the $(i,j)$-component 
\[P(G,K)_{i,j}=\Hom_{\Zz G}(P_i,C^j(\Rr\ZZ_K)),\ \ P_i=\Zz[m]_i\otimes\Zz H.\]
It is evident that $\PP(G,K)$ is functorial in a canonical way with respect to the pair $(G,K)$ for inclusions: $G\subset G'\subset H$, $K\subset K'\subset 2^{[m]}$. Similarly, for any coefficient ring $\kk$ we can define
\[\PP(G,K;\kk)=\Hom_{\kk G}(\PP(\kk),C^*(\Rr\ZZ_K;\kk)),\ \ \PP(\kk)=\Zz[m]\otimes\kk H.\]
\begin{prop}\label{prop:natral equivariant cohomology}
There is an isomorphism of graded $\kk$-modules for any subgroup $G\subset (\Zz/2)^m$ and any coefficient ring $\kk$: 
\[H^*(\PP(G,K;\kk))\cong H^*_G(\Rr\ZZ_K;\kk).\]
The isomorphism is natural with respect to the inclusions of subgroups of $(\Zz/2)^m$ and subcomplexes of $2^{[m]}$. 	
\end{prop}
\begin{proof}
	The proof is similar to the proof of Proposition \ref{prop:equivariant} and simpler for naturality by using $(S^\infty)^{\times m}$ as the total space of the universal $G$-bundle for any subgroup $G\subset(\Zz/2)^m$.
\end{proof}

Before proving Theorem \ref{thm:mod 2 image}, we also need several lemmas. Fix a coefficient ring $\kk$. We denote by $E_*^{*,*}(G,K)$ the spectral sequence of $\PP(G,K;\kk)$ corresponding to the decreasing filtration on $\PP(G,K;\kk)$ by the degree of $\PP(\kk)$, i.e., 
\[F_n\PP(G,K;\kk)=\bigoplus_{i\geq n}P(G,K;\kk)_{i,j}.\]
As usual, the first index $p$ in $E_*^{p,q}(G,K)$ indicates the filtration degree. $E_*^{*,*}(G,K)$ is functorial with respect to inclusions $G\subset G'$ and $K\subset K'$  in the sense of Proposition \ref{prop:natral equivariant cohomology}.
$E_*^{*,*}((\Zz/2)^m,K)$ is simply written as $E_*^{*,*}(K)$.
\begin{lem}\label{lem:E_infinity}
	The spectral sequence $E_*^{*,*}(K)$ for $\kk=\Zz/2$ satisfies $E_\infty^{*,*}(K)=E_\infty^{*,0}(K)=\Zz/2[K]$.
\end{lem}
\begin{proof}
	From \cite[Theorem 4.8]{DJ91}, we know that 
	\[H^*_{(\Zz/2)^m}(\Rr\ZZ_K;\Zz/2)\cong \Zz/2[K].\]
	In particular, we have $H^*_{(\Zz/2)^m}(\Rr\ZZ_{2^{[m]}};\Zz/2)\cong\Zz/2[m]$, and the inclusion $K\subset 2^{[m]}$ induces the  surjection $\Zz/2[m]\to\Zz/2[K]$ of mod $2$ $(\Zz/2)^m$-equivariant cohomology.
	
	By definition, the $E_1$-term of the spectral sequence $E_*^{*,*}(K)$ is
	\[E_1^{p,q}(K)=\Hom_{\Zz/2}(\Zz/2[m]_p,H^q(\Rr\ZZ_K;\Zz/2)).\]
	So $E_*^{*,*}(2^{[m]})$ collapses on the $E_1$-page:
	\[
	E_\infty^{*,*}(2^{[m]})=E_1^{*,*}(2^{[m]})=E_1^{*,0}(2^{[m]})=\Zz/2[m].\]
	Note that there are $\Zz/2$-module isomorphisms:
	\[
	E_\infty^{*,*}(K)\cong H^*(\PP((\Zz/2)^m,K;\Zz/2))\cong \Zz/2[K],\] 
	where the second isomorphism follows from Proposition \ref{prop:natral equivariant cohomology}. 
	Since $\Zz/2[m]\to\Zz/2[K]$ is surjective and $E_\infty^{*,*}(2^{[m]})=E_\infty^{*,0}(2^{[m]})$, a general fact in homological algebra (Lemma \ref{lem:homological algebra} below) shows that  $E_\infty^{*,*}(K)=E_\infty^{*,0}(K)=\Zz/2[K]$, and the proof is complete.
\end{proof}
\begin{lem}\label{lem:homological algebra}
Let $\DD=\oplus_{i,j\geq 0}D_{i,j}$ and $\DD'=\oplus_{i,j\geq 0}D'_{i,j}$ be first-quadrant double complexes with degree increasing differentials, and assume $\phi:\DD\to\DD'$ is a double complex map. Let $E_*^{*,*}(\DD)$ (resp. $E_*^{*,*}(\DD')$) be the spectral sequence associated to the filtration $F_n\DD=\oplus_{i\geq n}D_{i,j}$ (resp. $F_n\DD'=\oplus_{i\geq n}D'_{i,j}$).
If the induced cohomology map $\phi^*:H^*(\DD)\to H^*(\DD')$ is surjective and $E_\infty^{*,*}(\DD)=E_\infty^{*,0}(\DD)$, then $E_\infty^{*,*}(\DD')=E_\infty^{*,0}(\DD')$.
\end{lem}
\begin{proof}
By assumption, $H^n(\DD)=E_\infty^{n,0}(\DD)=F_nH^n(\DD)/F_{n+1}H^n(\DD)$. Since $\DD$ is a first-quadrant double complex, $F_{n+1}H^n(\DD)=0$, hence $H^n(\DD)=F_nH^n(\DD)$. So,
\[H^n(\DD')=\phi^*(H^n(\DD))=\phi^*(F_nH^n(\DD))\subset F_nH^n(\DD').\]
Thus, $H^n(\DD')=F_nH^n(\DD')=E_\infty^{n,0}(\DD')$ since $\DD'$ is also a first-quadrant double complex. This implies that $E_\infty^{*,*}(\DD')=E_\infty^{*,0}(\DD')$.
\end{proof}
We observe that $1\otimes u_\sigma\in\Zz/2[n]\otimes\RR_K=\DD(G,K)\otimes\Zz/2$ is a cocycle for any $\sigma\in K$ and any $G\cong(\Zz/2)^n\subset (\Zz/2)^m$. This is because $d(u_\sigma)=0$ in $\RR_K$ and $g\cdot u_\sigma\equiv u_\sigma$ mod $2$ for any $g\in (\Zz/2)^m$ by \eqref{eq:action}, so the differential of $\Zz/2[n]\otimes\RR_K$ maps $1\otimes u_\sigma$ to zero by \eqref{eq:differential D(G,K)}.
\begin{lem}\label{lem:homologous}
	For any face $\sigma\in K$, the cocycle $1\otimes u_\sigma$ is homologous to $ x_\sigma\otimes1$ in $\Zz/2[m]\otimes\RR_K$ with differentials given in \eqref{eq:differential D(G,K)} for the standard coordinate basis $\{g_1,\ldots,g_m\}$ of $(\Zz/2)^m$.
\end{lem}
\begin{proof}
	Suppose that $\sigma=\{i_1,\dots,i_k\}$, and let $\sigma_0=\emptyset$, $\sigma_j=\{i_1,\dots,i_j\}$, $\tau_j=\sigma\setminus\sigma_j$, $0\leq j\leq k$. Using \eqref{eq:action} and \eqref{eq:differential D(G,K)}, a direct calculation shows that for $1\leq j\leq k$
	\[d(x_{\tau_j}\otimes u_{\sigma_{j-1}}t_{i_j})\equiv x_{\tau_j}\otimes u_{\sigma_j} +x_{\tau_{j-1}}\otimes u_{\sigma_{j-1}} \bmod 2.\]
	Then we have
	\[d(\sum_{j=1}^k x_{\tau_j}\otimes u_{\sigma_{j-1}}t_{i_j})\equiv 1\otimes u_\sigma+x_\sigma\otimes1\bmod 2,\]
	and the lemma is proved.
  \end{proof}

\begin{proof}[Proof of Theorem \ref{thm:mod 2 image}]
	Because of the commutative diagram \eqref{diag:Phi and Psi}, we need only verify the statement for $\Phi_\Lambda$.
	Let $H=(\Zz/2)^m$ and $G=\ker\Lambda\cong(\Zz/2)^{m-n}$. Since $\Lambda:(\Zz/2)^m\to (\Zz/2)^n$ is a characteristic function over $K$, $G$ acts freely on $\Rr\ZZ_K$. So the projection \[EG\times_G\Rr\ZZ_K\to \Rr\ZZ_K/G=M(K,\Lambda)\]
	is a fiber bundle with fiber the contractible space $EG$, and then we have a homotopy equivalence $EG\times_G\Rr\ZZ_K\simeq M(K,\Lambda)$. Thus, for any coefficient ring $\kk$
	\[H^*(M(K,\Lambda);\kk)\cong H^*_G(\Rr\ZZ_K;\kk).\]
	Moreover, an argument similar to what we did in the first paragraph of the proof of Proposition \ref{prop:equivariant} shows that the cochain map
	\begin{equation}\label{eq:cochain map}
		\begin{split}
		C^*(\Rr\ZZ_K;\kk)^G=(\RR_K\otimes\kk)^G&\cong \Hom_{\kk G}(\kk,C^*(\Rr\ZZ_K;\kk))\\
		&\to \Hom_{\kk G}(\PP(\kk),C^*(\Rr\ZZ_K;\kk))=\PP(G,K;\kk)
	\end{split}
\end{equation}
	is a quasi-isomorphism that induces the above isomorphism.
	
	Since $K$ is Cohen-Macaulay over $\Zz/2$, there is the surjective map \[\Zz/2[K]\cong H^*_H(\Rr\ZZ_K;\Zz/2)\to H^*_{G}(\Rr\ZZ_K;\Zz/2)\cong H^*(M(K,\Lambda);\Zz/2)\]
    by Theorem \ref{thm:DJ}. Combining this with Proposition \ref{prop:natral equivariant cohomology} and Lemmas \ref{lem:E_infinity}, \ref{lem:homological algebra}, we get a commutative diagram
	\[\xymatrix{
		H^*_H(\Rr\ZZ_K;\Zz/2)\ar[r]^-{\cong}\ar[d]&E_\infty^{*,0}(H,K)\ar@{=}[r]\ar[d]&\Zz/2[K]\ar[d]\\
		H^*(M(K,\Lambda);\Zz/2)\ar[r]^-{\cong}&E_\infty^{*,0}(G,K)\ar@{=}[r]&\Zz/2[K]/\Theta_\Lambda}
	\]

	Recall that for any $\omega\in\row\Lambda$ and $\sigma\in K_\omega$, we have $\varphi_\omega(\sigma^*)\equiv u_\sigma$ mod $2$ by \eqref{eq:map phi_omega}.
	Since $g\cdot u_\sigma\equiv u_\sigma$ mod $2$ for any $g\in (\Zz/2)^m$, $u_\sigma$ is an element of both $(\RR_K\otimes\Zz/2)^H$ and 
	$(\RR_K\otimes\Zz/2)^G$. Consider another commutative diagram
	\[\xymatrix{
		(\RR_K\otimes\Zz/2)^H\ar^-{\Phi_H}[r]\ar[d]&\PP(H,K;\Zz/2)\cong\Zz/2[m]\otimes\RR_K\ar[d]\\
		(\RR_K\otimes\Zz/2)^G\ar^{\Phi_G}[r]&\PP(G,K;\Zz/2)}
	\]
	where the horizontal maps are given by \eqref{eq:cochain map}.
	Clearly, $\Phi_H(u_\sigma)=1\otimes u_\sigma$, which is homologous to $x_\sigma\otimes1$ in $\Zz/2[m]\otimes\RR_K$ by Lemma \ref{lem:homologous}.
	Hence $\Phi_H(u_\sigma)=x_\sigma$ in $H^*(\Zz/2[m]\otimes\RR_K)=\Zz/2[K]$. Since $\Phi_G$ is a quasi-isomorphism, the above two commutative diagrams show that $\varphi_\omega(\sigma^*)=x_\sigma$ in 
	\[H^*((\RR_K\otimes\Zz/2)^G)=H^*(M(K,\Lambda);\Zz/2)\cong \Zz/2[K]/\Theta_\Lambda.\] 
	Then the theorem follows from the fact that $\Phi_{\Lambda}$ is induced by $\bigoplus_{\omega\in\row\Lambda}\varphi_\omega$.
\end{proof}

\section{Real toric spaces associated to shellable complexes}\label{sec:shell}
A simplicial complex is said to be \emph{pure} if its facets (i.e. maximal faces)
all have the same dimension. We say that a pure simplicial complex $K$ is
\emph{shellable} if there is an ordering of its facets: $\sigma_1,\sigma_2,\dots,\sigma_s$ such that the following condition holds. Write $K_j$
for the subcomplex of $K$ generated by $\sigma_1,\dots,\sigma_j$, i.e.,
\[K_j=2^{\sigma_1}\cup2^{\sigma_2}\cup\cdots\cup2^{\sigma_j}.\]
Then we require that for all $1\leq i\leq s$, the set $K_i\setminus K_{i-1}$ of faces has a unique minimal element (with respect to inclusion), denoted $r(\sigma_i)$. (Formally, we set $K_0=\emptyset$, so $r(\sigma_1)=\emptyset$.) The ordering  $\sigma_1,\dots,\sigma_s$ of facets of $K$ is called a \emph{shelling} of $K$. It is well-known that a shellable simplicial complex
is Cohen-Macaulay over any field $\Ff$. 

The goal of this section is to prove the following result, which implies Theorem \ref{thm:group} when $K$ is a shellable sphere.

\begin{prop}\label{prop:shellable}
Let	$K$ be an $(n-1)$-dimensional shellable complex, and assume that $\Lambda:(\Zz/2)^m\to (\Zz/2)^n$ is a characteristic function over $K$. Then the maps $\Phi_\Lambda$, $\Psi_\Lambda$ in \eqref{eq:map} are surjective and are $2$-truncated isomorphisms (see below for the definition). 	
\end{prop}

Given an integer $k>0$ and a $\Zz$-module $M$, let
\[T_k(M)=\{a\in M\mid ka=0\}.\]
We say that a $\Zz$-module homomorphism $\phi:M_1\to M_2$ is a \emph{$k$-truncated isomorphism} if it induces an isomorphism of quotients
\[M_1/T_k(M_1)\cong M_2/T_k(M_2).\]
Since there is an obvious isomorphism $M/T_k(M)\xr{\cdot k}kM$ for any $\Zz$-module $M$, we have
\begin{lem}\label{lem:truncated}
For a $\Zz$-module homomorphism $\phi:M_1\to M_2$, the following are equivalent:
\begin{enumerate}[(a)]
	\item $\phi$ is a $k$-truncated isomorphism.
	\item $\phi:kM_1\to kM_2$ is an isomorphism.
\end{enumerate}
\end{lem}
The proof of Proposition \ref{prop:shellable} will use the following lemma that can be viewed as a generalization of one half of the five lemma.
\begin{lem}\label{lem:five lemma}
In a commutative diagram of $\Zz$-modules with exact rows,
\[\xymatrix{
	A\ar[r]^-{f}\ar[d]^-{\alpha}&B\ar[r]^-{g}\ar[d]^-{\beta}&C\ar[r]^-{h}\ar[d]^-{\gamma}&D\ar[r]^-{i}\ar[d]^-{\delta}&E\ar[d]^-{\varepsilon}\\
	A'\ar[r]^-{f'}&B'\ar[r]^-{g'}&C'\ar[r]^-{h'}&D'\ar[r]^-{i'}&E'}
\]
assume that the following conditions are satisfied
\begin{itemize}
	\item $\beta$ and $\delta$ are surjective and are $k$-truncated isomorphisms.
	\item $\alpha$ is surjective and $\varepsilon$ restricted to $kE$ is injective.
	\item $\im f\subset kB$ and $\im i\subset kE$.
	\item $g(b)\in kC$ implies $b\in kB$.
\end{itemize}
Then $\gamma$ is surjective and is a $k$-truncated isomorphism.
\end{lem}

\begin{proof}
	The proof is a standard diagram chase. First, we show that $\gamma$ is surjective. For any $c'\in C'$, $h'(c')=\delta(d)$ for some $d\in D$ since $\delta$ is surjective, and $i'h'(c')=0$ by exactness. Since $\im i\subset kE$ and $\varepsilon$ restricted to $kE$ is injective, the equation $\varepsilon i(d)=i'\delta(d)=i'h'(c')=0$ gives $i(d)=0$, hence $d=h(c)$ for some $c\in C$ by exactness of the upper row. Clearly, $h'$ maps the difference $c'-\gamma(c)$ to $0$.
	Therefore $c'-\gamma(c)=g'(b')$ for some $b'\in B'$ by exactness. Since $\beta$
	is surjective, $b'=\beta(b)$ for some $b\in B$, and then $\gamma(c+g(b))=\gamma(c)+g'\beta(b)=\gamma(c)+g'(b')=c'$, showing that $\gamma$ is surjective.
	
	To prove that $\gamma$ is a $k$-truncated isomorphism, it suffices to show that $\gamma$ restricted to $kC$ is injective by Lemma \ref{lem:truncated}, since we already know that it is surjective. Suppose that $\gamma(kc)=0$ for some $c\in C$. Since $\delta$ is a $k$-truncated isomorphism, $\delta h(kc)=h'\gamma(kc)=0$ implies $h(kc)=0$, so $kc=g(b)$ for some $b\in B$ by exactness, and then $b\in kB$ by the last condition in the lemma. The element $\beta(b)$ satisfies $g'\beta(b)=\gamma g(b)=\gamma(kc)=0$, so $\beta(b)=f'(a')$ for some $a'\in A'$.
	Since $\alpha$ is surjective, $a'=\alpha(a)$ for some $a\in A$. Since $\im f\subset kB$ and $b\in kB$, $b-f(a)\in kB$, and since $\beta$ is a $k$-truncated isomorphism, $\beta(b-f(a))=\beta(b)-\beta f(a)=\beta(b)-f'\alpha(a)=\beta(b)-f'(a')=0$ implies $b-f(a)=0$. Thus $b=f(a)$, and so $kc=g(b)=gf(a)=0$.
\end{proof}

Here we also give an analog of Lemma \ref{lem:five lemma}, which will be used in the next section.
\begin{lem}\label{lem:five lemma 2}
If the third and fourth conditions in Lemma \ref{lem:five lemma} are replaced by the following two conditions:
\begin{itemize}
	\item The maps $A'\otimes\Zz/k\xr{f'\otimes1} B'\otimes\Zz/k$ and $D'\otimes\Zz/k\xr{i'\otimes1} E'\otimes\Zz/k$ are injective.
	\item $\im h\subset kD$ and $T_k(C')\subset\im g'$.
\end{itemize}
Then the conclusion of Lemma \ref{lem:five lemma} still holds.
\end{lem}
\begin{proof}
	First, we show that $\gamma$ is surjective. For any $c'\in C'$, $i'h'(c')=0$ by exactness. The injectivity of the map $D'\otimes\Zz/k\xr{i'\otimes1} E'\otimes\Zz/k$ implies that $h'(c')=kd'$ for some $d'\in D'$. Since $\delta$ is surjective, there exists $d\in D$ such that $\delta(d)=d'$, so $\varepsilon i(kd)=i'(kd')=i'h'(c')=0$, and then $i(kd)=0$ since $\varepsilon$ restricted to $kE$ is injective. By exactness of the upper row, $kd=h(c)$ for some $c\in C$,  hence $h'$ maps the difference $c'-\gamma(c)$ to $0$. The rest is the same as in the first paragraph of the proof of Lemma \ref{lem:five lemma}.
	
	Next, we show that $\gamma$ restricted to $kC$ is injective. Suppose that $\gamma(kc)=0$ for some $c\in C$. Let $c'=\gamma(c)$. Then $c'\in T_k(C')$, hence $c'=g'(b')$ for some $b'\in B'$ by hypothesis. Since $\beta$ is surjective, $b'=\beta(b)$ for some $b\in B$, hence \[\gamma(c-g(b))=c'-\gamma(g(b))=c'-g'(b')=0,\]
	and so $\delta h(c-g(b))=0$. Since $\im h\subset kD$ and $\delta$ is a $k$-truncated isomorphism, we have $h(c-g(b))=0$ by Lemma \ref{lem:truncated}. Thus $c-g(b)=g(b_0)$ for some $b_0\in B$ by exactness, and so $c=g(b+b_0)$. Since
	 \[0=\gamma(kc)=\gamma g(k(b+b_0))=g'\beta(k(b+b_0)),\]
    $k\beta(b+b_0)=f'(a')$ for some $a'\in A'$ by exactness. Then the injectivity of the map  $A'\otimes\Zz/k\xr{f'\otimes1} B'\otimes\Zz/k$ implies that $a'=ka_0'$ for some $a_0'\in A'$. Since $\alpha$ is surjective, $a_0'=\alpha(a_0)$ for some $a_0\in A$, hence
	\[\beta(k(b+b_0-f(a_0)))=f'(a')-\beta f(ka_0)=f'(a')-f'(a')=0.\]
	It follows that $k(b+b_0-f(a_0))=0$ since $\beta$ is a $k$-truncated isomorphism. This, together with the fact that $c=g(b+b_0)$ and the exactness of the upper row, gives the desired equation $kc=0$. 
\end{proof}

The proof of Proposition \ref{prop:shellable} will also use the following fact.
\begin{lem}\label{lem:induction step}
Let	$K$ be an $(n-1)$-dimensional shellable complex with a shelling $\sigma_1,\dots,\sigma_s$, and let $K_j=\bigcup_{i=1}^j 2^{\sigma_i}$. Then for any $\omega\subset [m]$, we have
\[H^i((K_j)_\omega,(K_{j-1})_\omega)=
\begin{cases}
	\Zz\{r(\sigma_j)^*\},&\text{if }\omega\cap\sigma_j=r(\sigma_j)\text{ and } i=\dim r(\sigma_j),\\
	 0,&\text{otherwise.}
\end{cases}\]
Moreover, if $\Lambda:(\Zz/2)^m\to (\Zz/2)^n$ is a characteristic function over $K$, then for each $\sigma_j$ there exists a unique element $\omega_j\in\row\Lambda$ such that $\omega_j\cap\sigma_j=r(\sigma_j)$.
\end{lem}
\begin{proof}
	Since $K_j=K_{j-1}\cup 2^{\sigma_j}$, 
	\[H^*((K_j)_\omega,(K_{j-1})_\omega)\cong H^*(2^{\sigma_j\cap\omega},2^{\sigma_j}\cap (K_{j-1})_\omega)\] 
	by excision. Note that
	$2^{\sigma_j}\cap (K_{j-1})_\omega=\{\tau\subset\sigma_j\cap\omega\mid r(\sigma_j)\not\subset\tau\}$.
	Hence, if $\omega\cap\sigma_j=r(\sigma_j)$, then any proper subset of $r(\sigma_j)$ is a face of $(K_{j-1})_\omega$, so 
	\[H^*((K_j)_\omega,(K_{j-1})_\omega)=H^{\dim r(\sigma_j)}(2^{r(\sigma_j)},\partial r(\sigma_j))=\Zz\{r(\sigma_j)^*\}.\]
	 On the other hand, if $\omega\cap\sigma_j\neq r(\sigma_j)$, then it is easy to see that $2^{\sigma_j}\cap (K_{j-1})_\omega$ is a contractible complex. From here, the first assertion of the lemma follows easily.
	 
	 The final assertion is a consequence of \cite[Lemma 3.5]{CC21}.
	\end{proof}
	
	For the rest of this paper, we will use the abbreviated notations $M(K)$, $\A(K)$, $\RR(K)$, $\QQ(K)$ for $M(K,\Lambda)$, $\A(K,\Lambda)$, $\RR(K,\Lambda)$, $\QQ(K,\Lambda)$ respectively, whenever the characteristic function $\Lambda$ over $K$ is clear from the context. For a subcomplex $L\subset K$, let $\A(K,L)$, $\RR(K,L)$, $\QQ(K,L)$ be the kernels of the quotient maps 
	\[\A(K)\to\A(L),\ \ \RR(K)\to\RR(L),\ \ \QQ(K)\to\QQ(L),\]
	respectively. 

   We are now ready to prove our main result of this section. 
\begin{proof}[Proof of Proposition \ref{prop:shellable}]
	We will consider only the map $\Phi_\Lambda$ and the argument for $\Psi_\Lambda$ is very similar.
	The proof is by induction on the number of facets of $K$. Suppose that $\sigma_1,\dots,\sigma_s$ is a shelling of $K$, and let $K_j=\bigcup_{i=1}^j 2^{\sigma_i}$. 
	
	First consider the case $j=1$. By definition, $\Rr\ZZ_{K_1}=(D^1)^n\times (S^0)^{m-n}$ (we think of $K_1$ as a simplicial complex with $m-n$ ghost vertices), and \[M(K_1)=\Rr\ZZ_{K_1}/(\Zz/2)^{m-n}=D^n.\]
	So $H^*(\RR(K_1))=H^0(\RR(K_1))=\Zz$. Since $\Lambda$ is a characteristic function of rank $\dim K+1$, for an element $\omega\in\row\Lambda$, $ (K_1)_\omega=\emptyset$ only if $\omega=\emptyset$. It follows that 
	$H^*(\A(K_1))\cong\Zz\oplus A$,
	where the $\Zz$ summand corresponds to $\w H^*((K_1)_\emptyset)=\w H^{-1}(\emptyset)=\Zz\{\emptyset^*\}$ and $A$ is a direct sum of $(\Zz/2)$'s. With the notation in \eqref{eq:map phi_omega} and \eqref{eq:map}, we have $\Phi_\Lambda(\emptyset^*)=f_{\emptyset,\emptyset}=1$. Hence the theorem holds for $K_1$.
	
	Now assume that the theorem holds for $K_{j-1}$. 
	Consider the commutative diagram of long exact sequences
	\begin{equation}\label{diag:long exact}
		\begin{tikzcd}
		\cdots\rar{2\delta}& [-10pt]H^*(\A(K_j,K_{j-1}))\rar\dar{\Phi_\Lambda}&[-10pt]H^*(\A(K_j))\rar\dar{\Phi_\Lambda}&[-10pt]H^*(\A(K_{j-1}))\rar{2\delta}\dar{\Phi_\Lambda}&[-10pt]\cdots\\
		\cdots\rar{d}& H^*(\RR(K_j,K_{j-1}))\rar&H^*(\RR(K_j))\rar&H^*(\RR(K_{j-1}))\rar{d}&\cdots
	\end{tikzcd}
	\end{equation}
	By the induction hypothesis, the right-hand vertical map is surjective and is a $2$-truncated isomorphism. \cite[Proposition 2.5]{CC21} shows that up to homotopy, $M(K_j)$ is obtained from $M(K_{j-1})$ by attaching a cell of dimension $|r(\sigma_j)|$. So Theorem \ref{thm:isomorphism} gives
	\begin{equation}\label{eq:R(K_j,K_{j-1})}
		H^k(\RR(K_j,K_{j-1}))=
	\begin{cases}
		\Zz,&\text{if }k=|r(\sigma_j)|,\\
		0,&\text{otherwise.}
	\end{cases}
\end{equation}
	By Lemma \ref{lem:induction step} and the definition of $\A(K_j,K_{j-1})$, we have
	\[H^k(\A(K_j,K_{j-1}))/T_2(H^k(\A(K_j,K_{j-1})))=
	\begin{cases}
		\Zz\{r(\sigma_j)^*\},&\text{if } k=|r(\sigma_j)|,\\
		0,&\text{otherwise.}
	\end{cases}\]
Let $\omega_j\in\row\Lambda$ be the unique element such that $\omega_j\cap\sigma_j=r(\sigma_j)$.	It is easy to see that $0\neq f_{r(\sigma_j),\omega_j}\in \RR(K_j,K_{j-1})$ and that $\RR(K_j,K_{j-1})$ is zero in degrees less than $|r(\sigma_j)|$, so $f_{r(\sigma_j),\omega_j}$ is a nonzero multiple of a generator of $H^*(\RR(K_j,K_{j-1}))\cong \Zz$. However, since the coefficient of $u_{r(\sigma_j)}$ in the expression of $f_{r(\sigma_j),\omega_j}$ is $1$, $f_{r(\sigma_j),\omega_j}$ is actually a generator of $H^*(\RR(K_j,K_{j-1}))$. Thus, the left-hand vertical map is surjective and is a $2$-truncated isomorphism.

We now use Lemma \ref{lem:five lemma} to complete the induction step. The first and second conditions of Lemma \ref{lem:five lemma} have already been verified, and the third condition is automatically satisfied since the connecting homomorphism in the upper exact sequence of \eqref{diag:long exact} is $2\delta$. So it suffices to show that the map 
\[q^*:H^*(\A(K_j,K_{j-1}))\to H^*(\A(K_j))\]
satisfies the fourth condition in Lemma \ref{lem:five lemma}. Suppose that $q^*([b])=2[c]$ for some cohomology classes $[b]\in H^*(\A(K_j,K_{j-1}))$ and $[c]\in H^*(\A(K_j))$ with representatives $b\in\A(K_j,K_{j-1})$ and $c\in \A(K_j)$. Then $q^*(b)=2c+2\delta(a)$ for some $a\in \A(K_j)$, noting that the coboundary map of $\A(K_j)$ is $2\delta$. 
Hence $q^*(b)\in 2\A(K_j)$, and then $b\in 2\A(K_j,K_{j-1})$ since $\A(K_j,K_{j-1})$ is a direct summand of $\A(K_j)$. This shows that $[b]\in 2H^*(\A(K_j,K_{j-1}))$.
\end{proof}
We will need the following generalization of Proposition \ref{prop:shellable} in the next section.
\begin{prop}\label{prop:generalization} 
Let $L$ be a shellable complex, $K=2^{\sigma}*L$ for a simplex $2^\sigma$, and let $K'=\partial\sigma*L$. Assume that $\Lambda:(\Zz/2)^m\to (\Zz/2)^n$ ($n=\dim K+1$) is a characteristic function over $K$. Then the relative maps 
\begin{gather*}
	\Phi_\Lambda:H^*(\A(K,K'))\to H^*(\RR(K,K')),\\ \Psi_\Lambda:H^*(\A(K,K'))\to H^*(\QQ(K,K'))
\end{gather*} 
are surjective and are $2$-truncated isomorphisms. 
\end{prop}
\begin{proof}
The proof is similar to the proof of Proposition \ref{prop:shellable}, using induction on the number of facets of $L$, and we consider only the map $\Phi_\Lambda$. Suppose that $\tau_1,\dots,\tau_s$ is a shelling of $L$, and let $L_j=\bigcup_{i=1}^j 2^{\tau_i}$, $K_j=2^{\sigma}*L_j$, $K_j'=\partial\sigma*L_j$ for $1\leq j\leq s$. The base case $(K_1,K_1')$ can be proved in the same way as  for the left-hand vertical map in \eqref{diag:long exact}. For the induction step, let 
\[K_j''=(\partial\sigma*L_j)\cup(2^\sigma*L_{j-1})=K_j'\cup K_{j-1}\]
and consider the commutative diagram of long exact sequences
\[\begin{tikzcd}
	\cdots\rar{2\delta}& [-10pt]H^*(\A(K_j,K_j''))\rar\dar{\Phi_\Lambda}
	&[-10pt]H^*(\A(K_j,K_j'))\rar\dar{\Phi_\Lambda}
	&[-10pt]H^*(\A(K_j'',K_j'))\rar{2\delta}\dar{\Phi_\Lambda}&[-10pt]\cdots\\
	\cdots\rar{d}&H^*(\RR(K_j,K_j''))\rar&H^*(\RR(K_j,K_j'))\rar&H^*(\RR(K_j'',K_j'))\rar{d}&\cdots
\end{tikzcd}
\]
By definition, the pair $(K_j'',K_j')$ can be replaced by $(K_{j-1},K_{j-1}')$.
Hence the right-hand vertical map in the above diagram is surjective and is a $2$-truncated isomorphism by the induction hypothesis. Set $\rho_j=\sigma\cup r(\tau_j)$. Then $\rho_j$ is the unique minimal element in $K_j\setminus K_j''$, which implies that the results of Lemma \ref{lem:induction step} and \eqref{eq:R(K_j,K_{j-1})} also hold for the pair $(K_j,K_j'')$ with $r(\sigma_j)$ replaced by $\rho_j$. Hence, arguing as in the proof of Proposition \ref{prop:shellable}, we can prove that the left-hand vertical map is surjective and is a $2$-truncated isomorphism.

Finally, using Lemma \ref{lem:five lemma}, one can show that the middle vertical map in the above diagram is surjective and is a $2$-truncated isomorphism (see the last paragraph of the proof of Proposition \ref{prop:shellable}), finishing the induction step.
\end{proof}
\section{Polytopal spheres and subdivisions of PL spheres}\label{sec:PL}
A simplicial sphere is said to be \emph{polytopal} if it is combinatorially equivalent to the boundary complex of a simplicial polytope.
We have the following inclusion relations between subclasses of simplicial spheres:
\[\text{polytopal spheres}\subset \text{shellable spheres}\subset \text{PL spheres},\]
where the first inclusion is a famous result of Bruggesser and Mani \cite{BM72}. Moreover, in dimensions greater than $2$, these inclusions are strict. 

The following combinatorial result of Adiprasito and Izmestiev plays a crucial role in the proof of Theorem \ref{thm:group}.
\begin{thm}[{\cite[Theorem 1]{AI15}}]\label{thm:subdivision}
	Let $K$ be a PL sphere. Then there exists a polytopal sphere $K'$
	such that $K'$ is obtained from $K$ by a sequence of stellar subdivisions.
\end{thm}
The proof of Theorem \ref{thm:group} is by induction on the dimension of PL spheres, and we only deal with the map $\Phi_\Lambda$ as before. We first establish several preliminary results.

Let $K$ be a simplicial $(n-1)$-sphere on $[m]$, $L\subset K$ be a simplicial  $(n-1)$-ball. For each $i\geq 0$, let $\Delta^{i-1}$ be the $(i-1)$-simplex ($\Delta^{-1}=\{\emptyset\}$) with vertices $m+1,\dots,m+i$, and define 
\begin{equation}\label{eq:K^i}
\begin{gathered}
K^i=\Delta^{i-1}*K,\quad\partial K^i=\partial \Delta^{i-1}*K,\\
\Gamma^i=\partial K^i\cup(\Delta^{i-1}*L).
\end{gathered}
\end{equation}
Assume that $\Lambda:(\Zz/2)^{m+i}\to (\Zz/2)^{n+i}$ is a characteristic function over $K^i$.
Then we have a commutative diagram of long exact sequences
\[\begin{tikzcd}
	\cdots\rar{2\delta}& [-10pt]H^*(\A(K^i,\Gamma^i))\rar\dar{\Phi_\Lambda}
	&[-10pt]H^*(\A(K^i,\partial K^i))\rar\dar{\Phi_\Lambda}
	&[-10pt]H^*(\A(\Gamma^i,\partial K^i))\rar{2\delta}\dar{\Phi_\Lambda}&[-10pt]\cdots\\
	\cdots\rar{d}&H^*(\RR(K^i,\Gamma^i))\rar&H^*(\RR(K^i,\partial K^i))\rar&H^*(\RR(\Gamma^i,\partial K^i))\rar{d}&\cdots
\end{tikzcd}
\]

\begin{prop}\label{prop:iff}
	In the above diagram, if the left-hand vertical map is surjective and is a $2$-truncated isomorphism for all $*\geq 0$, then the middle vertical map is surjective and is a $2$-truncated isomorphism for all $*\geq 0$ if and only if the same is true for the right-hand vertical map.
\end{prop}
The `if' part of this proposition is a consequence of Lemma \ref{lem:five lemma}, as we have seen in the proof of Proposition \ref{prop:shellable}. For the `only if' part  we will need some fundamental results in algebraic combinatorics.

Let $K$ be an $(n-1)$-dimensional simplicial complex. Then there exists a sequence of algebraically independent homogeneous polynomials $\theta_1,\dots,\theta_n\in\Zz/2[K]$ such that $\Zz/2[K]$ is a
finitely generated $\Zz/2[\theta_1,\dots,\theta_n]$-module. Such a sequence $\theta_1,\dots,\theta_n$ is called a \emph{homogeneous system of parameters} (briefly h.s.o.p.) of $\Zz/2[K]$. If all $\theta_i$ are linear, then it is called a \emph{linear system of parameters} (briefly l.s.o.p.). An l.s.o.p. of $\Zz/2[K]$ may fail to exist for some simplicial complex $K$ (see \cite[Example 3.3.4]{BP15}). Here is a useful criterion for l.s.o.p.
\begin{prop}[see {\cite[Theorem 5.1.16]{BH93}}]\label{prop:lsop}
Let $K$ be an $(n-1)$-dimensional simplicial complex. Assume that $\Lambda:(\Zz/2)^m\to (\Zz/2)^n$ is a surjective linear map, and denote by $\lambda_i$ the $i$th row vector of the matrix corresponding to $\Lambda$. Then the following are equivalent:
	\begin{enumerate}[(a)]
		\item $\Lambda$ is a characteristic function over $K$.
		\item The sequence of linear forms 
		$\theta_i=\sum_{j\in\lambda_i}x_j$, $1\leq i\leq n$,
		is an l.s.o.p. of $\Zz/2[K]$.
	\end{enumerate}
\end{prop}
Let $\theta_1,\dots,\theta_n$ be an l.s.o.p. of $\Zz/2[K]$, and let $\Theta$ be the ideal of $\Zz/2[K]$ generated by $\theta_i$'s. The quotient ring $\Zz/2[K]/\Theta$ is called an \emph{Artinian reduction} of $\Zz/2[K]$ and is denoted $\Zz/2(K;\Theta)$ or simply $\Zz/2(K)$ when there is no risk of confusion. For a subcomplex $L\subset K$, let $I(K,L)$ be the ideal of $\Zz/2[K]$ generated by monomials in $\{x_\sigma\mid \sigma\in K\setminus L\}$. Then there is a short exact sequence
\[0\to I(K,L)\to \Zz/2[K]\to \Zz/2[L]\to 0.\]
If we quotient out by the ideal $\Theta\subset \Zz/2[K]$ we obtain another short exact sequence
\begin{equation}\label{eq:exact}
	0\to I(K,L)/(I(K,L)\cap\Theta)\to \Zz/2(K;\Theta)\to \Zz/2[L]/\Theta\to 0.
\end{equation}
 Since $I(K,L)\Theta\subset I(K,L)\cap\Theta$, there is a surjection
\[I(K,L)/I(K,L)\Theta\to I(K,L)/(I(K,L)\cap\Theta).\]
Let $\Zz/2(K,L;\Theta)$ (or simply $\Zz/2(K,L)$) denote the quotient $I(K,L)/I(K,L)\Theta$. The following result shows that the above surjection becomes an equality in some cases. 

\begin{prop}\label{prop:short exact}
	Let $K$ be a simplicial sphere and $L\subset K$ be a simplicial ball of the same dimension. Then for any Artinian reduction $\Zz/2(K)$, there is a short exact sequence 
	\[0\to\Zz/2(K,L)\to\Zz/2(K)\to\Zz/2(L)\to0.\]
(Note that the image of an l.s.o.p. of $\Zz/2[K]$ in $\Zz/2[L]$ is also an l.s.o.p. of $\Zz/2[L]$ by Proposition \ref{prop:lsop}.)
\end{prop}
\begin{proof}
	See the first paragraph of the proof of \cite[Theorem 3.1]{Swa14}.
\end{proof}
\begin{lem}\label{lem:relative isomorphism}
	Let $K^i$, $\Gamma^i$ be as in \eqref{eq:K^i}. Then for any characteristic function $\Lambda:(\Zz/2)^{m+i}\to (\Zz/2)^{n+i}$ over $K^i$, there are isomorphisms:
	\begin{gather*}
		H^*(\RR(K^i,\Gamma^i)\otimes\Zz/2)\cong \Zz/2(K^i,\Gamma^i),\\ H^*(\RR(K^i,\partial K^i)\otimes\Zz/2)\cong \Zz/2(K^i,\partial K^i),
	\end{gather*}
	and the map $\Zz/2(K^i,\Gamma^i)\to \Zz/2(K^i,\partial K^i)$ is injective. (The ideal $\Theta_\Lambda\subset\Zz/2[K^i]$ generated by the row vectors of $\Lambda$ is omitted from the notation.) 
\end{lem}
\begin{proof}
	Let $\wh K^i=\partial\Delta^i*K$ be the $(n+i-1)$-dimensional simplicial sphere on $[m+i+1]$ obtained from $K^i$ by attaching the $(n+i-1)$-dimensional simplicial ball $B^i:=2^{\{m+i+1\}}*\partial K^i$ along $\partial K^i$. Let $\Lambda':(\Zz/2)^{m+i+1}\to (\Zz/2)^{n+i}$ be the linear map such that $\Lambda'(\bm a_j)=\Lambda(\bm a_j)$ for $1\leq j\leq m+i$, and
	\[\Lambda'(\bm a_{m+i+1})=\sum_{j=m+1}^{m+i}\Lambda(\bm a_j),\] 
	where $\bm a_j$ is the $j$th coordinate vector of $(\Zz/2)^{m+i+1}$. It is easy to check that $\Lambda'$ is a characteristic function  over $\wh K^i$. 
	Let $\Theta_{\Lambda'}\subset\Zz/2[\wh K^i]$ be the ideal generated by the row vectors of $\Lambda'$. 
	
	Define $\wh\Gamma^i=B^i\cup \Gamma^i$, which is an $(n+i-1)$-dimensional simplicial ball in $\wh K^i$. One easily sees that
	\begin{equation}\label{eq:relative isomorphism}
		\Zz/2(K^i,\Gamma^i)\cong \Zz/2(\wh K^i,\wh \Gamma^i)\ \text{ and }\ \Zz/2(K^i,\partial K^i)\cong \Zz/2(\wh K^i,B^i),
	\end{equation}
	where the notation $\Theta_{\Lambda'}$ is omitted from the right side of the isomorphisms.
	Since $\wh K^i$ and $\wh\Gamma^i$ are Cohen-Macaulay over $\Zz/2$,  we have isomorphisms:
	\[
		H^*(M(\wh K^i,\Lambda');\Zz/2)\cong \Zz/2(\wh K^i)\ \text{ and }\  H^*(M(\wh\Gamma^i,\Lambda');\Zz/2)\cong \Zz/2(\wh\Gamma^i)
	\]
	by Theorem \ref{thm:DJ}. Hence the natural map $H^*(M(\wh K^i,\Lambda');\Zz/2)\to H^*(M(\wh\Gamma^i,\Lambda');\Zz/2)$ is surjective, and
	the mod $2$ relative cohomology for the pair 
	$(M(\wh K^i,\Lambda'),M(\wh\Gamma^i,\Lambda'))$ is isomorphic to the kernel of the map $\Zz/2(\wh K^i)\to \Zz/2(\wh\Gamma^i)$. Combining this with Proposition \ref{prop:short exact} and the first isomorphism in \eqref{eq:relative isomorphism} gives 
	\[H^*(M(\wh K^i,\Lambda'),M(\wh\Gamma^i,\Lambda');\Zz/2)\cong \Zz/2(K^i,\Gamma^i).\]
	Note that
	\[\begin{split}
		H^*(\RR(K^i,\Gamma^i)\otimes\Zz/2)&\cong H^*(M(K^i,\Lambda),M(\Gamma^i,\Lambda);\Zz/2)\ \text{ by Theorem \ref{thm:isomorphism}}\\
		&\cong H^*(M(K^i,\Lambda'),M(\Gamma^i,\Lambda');\Zz/2)\\
		&\cong H^*(M(\wh K^i,\Lambda'),M(\wh\Gamma^i,\Lambda');\Zz/2)\ \text{ by excision,}
	\end{split}\]
	where the second isomorphism comes from the homeomorphisms of toric spaces: \[M(K^i,\Lambda)\cong M(K^i,\Lambda')\ \text{ and }\ M(\Gamma^i,\Lambda)\cong M(\Gamma^i,\Lambda'),\] 
	viewing the right side as toric spaces over $K^i$ and $\Gamma^i$ respectively with $\{m+i+1\}$ as a ghost vertex.
	Hence, the first isomorphism in the lemma holds. The second isomorphism can be proved in the same way, using the pair $(M(\wh K^i,\Lambda'),M(B^i,\Lambda'))$.
	
	For the second assertion of the lemma, by \eqref{eq:relative isomorphism} it suffices to show that the map $\Zz/2(\wh K^i,\wh\Gamma^i)\to \Zz/2(\wh K^i,B^i)$ is injective. By \eqref{eq:exact} and Proposition \ref{prop:short exact}, we have
	\begin{gather*}
		\Zz/2(\wh K^i,\wh\Gamma^i)\cong I(\wh K^i,\wh\Gamma^i)/(I(\wh K^i,\wh\Gamma^i)\cap\Theta_{\Lambda'}),\\ 
		\Zz/2(\wh K^i,B^i)\cong I(\wh K^i,B^i)/(I(\wh K^i,B^i)\cap\Theta_{\Lambda'}).
	\end{gather*}
	Since $B^i\subset \wh\Gamma^i$, there is an inclusion
	\[I(\wh K^i,\wh\Gamma^i)/(I(\wh K^i,\wh\Gamma^i)\cap\Theta_{\Lambda'})\subset I(\wh K^i,B^i)/(I(\wh K^i,B^i)\cap\Theta_{\Lambda'}),\]
which gives the desired conclusion.
	\end{proof}
We are now ready to prove the `only if' part of Proposition \ref{prop:iff}.
\begin{proof}[Proof of Proposition \ref{prop:iff} ($\Rightarrow$)]
For the triple $(K^i,\Gamma^i,\partial K^i)$, consider the cohomology maps:
\[H^*(\RR(K^i,\Gamma^i))\xr{f^*}H^*(\RR(K^i,\partial K^i))\xr{g^*}H^*(\RR(\Gamma^i,\partial K^i)).\] 
Viewing the map 
$\Phi_\Lambda:H^*(\A(\Gamma^i,\partial K^i))\to H^*(\RR(\Gamma^i,\partial K^i))$ as $\gamma$ in Lemma \ref{lem:five lemma 2}, we only need to verify that the map
\[f^*\otimes 1:H^*(\RR(K^i,\Gamma^i))\otimes\Zz/2\to H^*(\RR(K^i,\partial K^i))\otimes\Zz/2\]
is injective and that $T_2(H^*(\RR(\Gamma^i,\partial K^i)))\subset\im g^*$, since the other conditions in Lemma \ref{lem:five lemma 2} automatically hold by assumption.
Indeed, the injectivity of $f^*\otimes 1$ follows easily from Lemma \ref{lem:relative isomorphism} and the universal coefficient theorem.

Since the mod $2$ cohomology map \[f^*:H^*(\RR(K^i,\Gamma^i)\otimes\Zz/2)\to H^*(\RR(K^i,\partial K^i)\otimes\Zz/2)\] is injective by Lemma \ref{lem:relative isomorphism}, the lower long exact sequence in the diagram preceding Proposition \ref{prop:iff} shows that the map 
\[g^*:H^*(\RR(K^i,\partial K^i)\otimes\Zz/2)\to H^*(\RR(\Gamma^i,\partial K^i)\otimes\Zz/2)\]
is surjective. Hence the following lemma gives the desired inclusion \[T_2(H^*(\RR(\Gamma^i,\partial K^i)))\subset\im g^*.\] 
\end{proof}
\begin{lem}
	Let $B$ and $C$ be two chain complexes of free abelian groups with homology of finite type and suppose that there is a chain map $g:B\to C$. If, for a prime $p$, the induced homology map 
	\[g_*:H_*(B\otimes\Zz/p)\to H_*(C\otimes\Zz/p)\]
	is surjective, then $T_p(H_*(C))\subset\im g_*$.
\end{lem}
\begin{proof}
	Let $[c]\in T_p(H_k(C))$ be a nonzero homology class and let $l\geq 0$ be the maximal integer such that $[c]\in p^lH_k(C)$. Write $[c]=p^l[c_0]$, so that $p^{l+1}[c_0]=0$ and $[c_0]\not\equiv 0$ mod $p$ in $H_k(C)$. Then there exists $c_1\in C_{k+1}$ such that $dc_1=p^{l+1}c_0$. Indeed, $0\neq[c_1]\in H_{k+1}(C\otimes\Zz/p)$ since $d_{l+1}([c_1])=[c_0]\neq 0$ in the $E_{l+1}$-term of the Bockstein spectral sequence of $H_*(C\otimes\Zz/p)$.
	The condition on $g_*$ implies that there is $b_1\in B_{k+1}$ such that $db_1\in pB$ and $g(b_1)=c_1+pc_1'+dc_2$ for some $c_1'\in C_{k+1}$ and $c_2\in C_{k+2}$. Set $b=\frac{db_1}{p}$. Then 
	\[g(b)=\frac{dg(b_1)}{p}=\frac{dc_1+pdc_1'}{p}=p^lc_0+dc_1'.\] 
	Thus, we have  $g_*([b])=[c]$. This completes the proof.
\end{proof}
Finally, we prove Theorem \ref{thm:group} by induction on the dimension of PL spheres, using the following proposition. The induction starts with the trivial case of the $(-1)$-dimensional sphere, i.e., $K=\{\emptyset\}$.
\begin{prop}\label{prop:induction}
	In the notation of \eqref{eq:K^i}, if for all PL spheres $K$ of dimension less than $n$ the map 
	\[\Phi_\Lambda:H^*(\A(K^i,\partial K^i))\to H^*(\RR(K^i,\partial K^i))\]
	is surjective and is a $2$-truncated isomorphism for any $i\geq 0$ and any characteristic function over $K^i$, then the same is true for all $n$-dimensional PL spheres.
\end{prop}
\begin{proof}
Suppose that $K$ is a PL $n$-sphere. By Theorem \ref{thm:subdivision}, there is a sequence of PL $n$-spheres: 
$K=K_0,K_1,\dots,K_s$
such that $K_s$ is polytopal  
and for $0\leq j\leq s-1$, $K_{j+1}=S_{\sigma_j}K_j$ for some nonempty face $\sigma_j\in K_j$. 

For each $i\geq 0$ and $0\leq j\leq s$, let 
\[K_j^i=\Delta^{i-1}*K_j,\ \ \partial K_j^i=\partial\Delta^{i-1}*K_j.\] 
Assume that $\Lambda_0:(\Zz/2)^m\to (\Zz/2)^{n+i+1}$ is an arbitrary characteristic function over $K^i_0=K^i$, with the convention that the vertex set of $K^i$ is $[m]$ for all $i\geq 0$. Number the new vertices $v_{\sigma_0},\dots,v_{\sigma_{s-1}}$  as $m+1,\dots,m+s$, and inductively, let $\Lambda_j:(\Zz/2)^{m+j}\to (\Zz/2)^{n+i+1}$ be the linear map such that $\Lambda_j(\bm a_k)=\Lambda_{j-1}(\bm a_k)$ for $1\leq k< m+j$, and
\[\Lambda_j(\bm a_{m+j})=\sum_{k\in\sigma_{j-1}}\Lambda_{j-1}(\bm a_k),\] 
where $\bm a_k$ is the $k$th coordinate vector of $(\Zz/2)^{m+j}$. Then $\Lambda_j$ is a characteristic function over $K_j^i$. For the rest of the proof, $\Lambda_j$ will be omitted from the notation.

Since polytopal spheres are shellable, the map 
\[\Phi_\Lambda:H^*(\A(K^i_s,\partial K^i_s))\to H^*(\RR(K^i_s,\partial K^i_s))\] 
is surjective and is a $2$-truncated isomorphism for $i\geq 0$ by Proposition \ref{prop:generalization}. Assume inductively that this is true for the pair $(K^i_{j+1},\partial K^i_{j+1})$. Let $L_j$ be the closure of $K_j\setminus \st_{K_j}\sigma_j$, which is a simplicial $n$-ball in $K_j$. It is easy to see that $L_j$ is also the closure of $K_{j+1}\setminus \st_{K_{j+1}}\{v_{\sigma_j}\}$
by the definition of stellar subdivision. Let
\[\Gamma^i_j=(\Delta^{i-1}*L_j)\cup \partial K_j^i\subset K_j^i,\]
and consider the commutative diagram of long exact sequences
\[\begin{tikzcd}
	\cdots\rar{2\delta}& [-10pt]H^*(\A(K^i_j,\Gamma^i_j))\rar\dar{\Phi_\Lambda}
	&[-12pt]H^*(\A(K^i_j,\partial K^i_j))\rar\dar{\Phi_\Lambda}
	&[-12pt]H^*(\A(\Gamma^i_j,\partial K^i_j))\rar{2\delta}\dar{\Phi_\Lambda}&[-12pt]\cdots\\
	\cdots\rar{d}& H^*(\RR(K^i_j,\Gamma^i_j))\rar&H^*(\RR(K^i_j,\partial K^i_j))\rar&H^*(\RR(\Gamma^i_j,\partial K^i_j))\rar{d}&\cdots
\end{tikzcd}
\]
We need to prove that the middle vertical map in the above diagram is surjective and is a $2$-truncated isomorphism.

Let $B_j^i$ be the $(n+i)$-ball $\Delta^{i-1}*\st_{K_j}\sigma_j$ with boundary $\partial B_j^i$. Then we have
\[\A(K^i_j,\Gamma^i_j)\cong \A(B_j^i,\partial B_j^i)\ \text{ and }\ \RR(K^i_j,\Gamma^i_j)\cong\RR(B_j^i,\partial B_j^i).\] 
Note that
\[(B_j^i,\partial B_j^i)\cong (\Delta^{i+|\sigma_j|-1}*\lk_{K_j}\sigma_j,\,\partial\Delta^{i+|\sigma_j|-1}*\lk_{K_j}\sigma_j).\]
Since the dimension of $\lk_{K_j}\sigma_j$ is less than $n$, the map $\Phi_\Lambda$ associated to the pair $(B_j^i,\partial B_j^i)$ is surjective and is a $2$-truncated isomorphism by assumption. Hence, to complete the proof, it suffices to show that the same is true for the right-hand vertical map in the above diagram by Proposition \ref{prop:iff}.

Since $L_j\subset K_{j+1}$, $D_j^i:=(\Delta^{i-1}*L_j)\cup\partial K_{j+1}^i$ is a subcomplex of $K_{j+1}^i$. Set 
\[C_j^i=\Delta^{i-1}*\st_{K_{j+1}}\{v_{\sigma_j}\}= \Delta^i*\lk_{K_{j+1}}\{v_{\sigma_j}\},\]
which is an $(n+i)$-ball in $K_{j+1}^i$, and let $\partial C_j^i= \partial\Delta^i*\lk_{K_{j+1}}\{v_{\sigma_j}\}$ be the boundary sphere of $C_j^i$.
By definition, one easily checks that 
\begin{gather}
\A(K_{j+1}^i,D_j^i)\cong \A(C_j^i,\partial C_j^i),\quad \RR(K_{j+1}^i,D_j^i)\cong \RR(C_j^i,\partial C_j^i),\label{eq:excision1}\\
\A(D_j^i,\partial K_{j+1}^i)\cong \A(\Gamma^i_j,\partial K^i_j),\quad \RR(D_j^i,\partial K_{j+1}^i)\cong \RR(\Gamma^i_j,\partial K^i_j).\label{eq:excision2}
\end{gather}
Hence, if we replace $K^i_j$, $\Gamma^i_j$, $\partial K_j^i$ with $K_{j+1}^i$, $D_j^i$, $\partial K_{j+1}^i$, respectively, in the above diagram, then the left-hand vertical map is surjective and is a $2$-truncated isomorphism by assumption and by \eqref{eq:excision1} since $\lk_{K_{j+1}}\{v_{\sigma_j}\}$ has dimension $n-1$. The same is true for the middle vertical map by induction. Then, by the `only if' part of Proposition \ref{prop:iff} and by \eqref{eq:excision2}, we get the desired result for the map $\Phi_\Lambda$ associated to the pair $(\Gamma^i_j,\partial K^i_j)$. 
\end{proof}

\section{The proofs of Theorems \ref{thm:ring} and \ref{thm:ring isomorphism}}\label{sec:ring}
In this section, $K$ is an $(n-1)$-dimensional simplicial complex and $\Lambda:(\Zz/2)^m\to (\Zz/2)^n$ is a characteristic function over $K$. The following two lemmas will be used to prove Theorems \ref{thm:ring} and \ref{thm:ring isomorphism}.

\begin{lem}\label{lem:product of h}
	Given $h_{\sigma_1,\omega_1}$, $h_{\sigma_2,\omega_2}\in\QQ(K,\Lambda)$ with $\omega_1,\omega_2\in\row\Lambda$, such that $\sigma_1\cap\omega_2=\emptyset$ and $\sigma_2\cap\omega_1=\emptyset$, their product in $\QQ(K,\Lambda)$ is given by
	\[h_{\sigma_1,\omega_1} h_{\sigma_2,\omega_2}=
		h_{\sigma_1\cup\sigma_2,\omega_1+\omega_2}.\]
	Here $\sigma_1\cup\sigma_2$ means the juxtaposition of the two ordered faces $\sigma_1$ and $\sigma_2$.
\end{lem}

\begin{lem}\label{lem:coboundary}
	Given $\omega_1,\omega_2\in\row\Lambda$ and two cocycles 
	\[c_i=\sum_{\sigma\in K_{\omega_i}}l_{i,\sigma}\sigma^*\in \w C^*(K_{\omega_i}),\ \ l_{i,\sigma}\in\Zz,\ \ i=1,2,\] 
	let
	\[
	c_1'=\sum_{\sigma\in K_{\omega_1},\,\sigma\cap\omega_2=\emptyset}l_{1,\sigma}\sigma^*,\quad c_2'=\sum_{\sigma\in K_{\omega_2},\,\sigma\cap\omega_1=\emptyset}l_{2,\sigma}\sigma^*,\]
	and set \[\alpha=\psi_{\omega_1}(c_1)\psi_{\omega_2}(c_2)-\psi_{\omega_1}(c_1')\psi_{\omega_2}(c_2').\]
 Then $2\alpha$ is a coboundary in $\QQ(K,\Lambda)$.
\end{lem}
\begin{proof}[Proof of Theorem \ref{thm:ring}]
With the notation of Lemma \ref{lem:coboundary}, let $\omega_1'=\omega_1\setminus\omega_2$ and $\omega_2'=\omega_2\setminus\omega_1$, and let $i_j:K_{\omega_j'}\hookrightarrow K_{\omega_j}$, $j=1,2$, be the simplicial inclusions. It is clear that $i_j^*(c_j)=c_j'$, $j=1,2$. By Lemma \ref{lem:product of h}, $\psi_{\omega_1}(c_1')\psi_{\omega_2}(c_2')=\psi_{\omega_1+\omega_2}(c_1'*c_2')$, where `$*$' denotes the product induced by the simplicial inclusion $K_{\omega_1+\omega_2}\to K_{\omega_1'}*K_{\omega_2'}$. Lemma \ref{lem:coboundary} implies that $\psi_{\omega_1}(c_1)\psi_{\omega_2}(c_2)=\psi_{\omega_1}(c_1')\psi_{\omega_2}(c_2')$ in $H^*(\QQ(K,\Lambda))/\QQ_2$. Combining these facts, we see that the isomorphism 
\[\Psi_\Lambda:H^*(\A(K,\Lambda))/\A_2\cong H^*(\QQ(K,\Lambda))/\QQ_2\]
in Theorem \ref{thm:group} is actually a ring isomorphism when the group \[H^*(\A(K,\Lambda))/\A_2\cong \GG^*(K,\Lambda)=\bigoplus_{\omega\in \row\Lambda}\w H^{*-1}(K_\omega)\]
is equipped with the $*$-product, where the group isomorphism is implied by the fact that the differential of $\A^*(K_\omega)$ is twice the differential of $\w C^*(K_\omega)$. This completes the proof of Theorem \ref{thm:ring}, using Theorem \ref{thm:DGA}. 
\end{proof}
Finally, we conclude this section with proofs of the preceding two lemmas.
\begin{proof}[Proof of Lemma \ref{lem:product of h}]
	If $\sigma_1\cap\omega_2=\emptyset$ and $\sigma_2\cap\omega_1=\emptyset$, then by the definition of  $h_{\sigma,\omega}$ (see \eqref{eq:f,h}), the product of $h_{\sigma_1,\omega_1}$ and $h_{\sigma_2,\omega_2}$ takes the form
	\[h_{\sigma_1,\omega_1} h_{\sigma_2,\omega_2}=\sum_{\substack{\omega\subset\omega_1\cup\omega_2,\\(\sigma_1\cup\sigma_2)\cap\omega=\emptyset}}l_\omega z_{\sigma_1\cup\sigma_2}y_{\omega},\quad l_\omega\in\Zz.\]
	We need to show that $l_\omega=0$ if $\omega\cap(\omega_1\cap\omega_2)\neq\emptyset$, and $l_\omega=(-2)^{|\omega|}$ otherwise. 
	
	Given $\omega\subset\omega_1\cup\omega_2$ with $(\sigma_1\cup\sigma_2)\cap\omega=\emptyset$, assume that $|\omega\cap(\omega_1\cap\omega_2)|=k$. Set
	\[P_i=\{(\omega_1',\omega_2')\mid \omega_1'\subset\omega_1,\,\omega_2'\subset\omega_2,\, \omega_1'\cup\omega_2'=\omega,\,|\omega_1'\cap\omega_2'|=i\},\ 0\leq i\leq k.\]
	Then $|P_i|=2^{k-i}\binom{k}{i}$, where the binomial coefficient $\binom{k}{i}$ is the number of choices of $\omega_1'\cap\omega_2'$ from $\omega\cap(\omega_1\cap\omega_2)$ and the number $2^{k-i}$ corresponds to the complement of $(\omega_1'\cap\omega_2')$ in $\omega\cap(\omega_1\cap\omega_2)$ with each element either in $\omega_1'$ or in $\omega_2'$. Each $(\omega_1',\omega_2')\in P_i$ gives a pair $(z_{\sigma_1}y_{\omega_1'},z_{\sigma_2}y_{\omega_2'})$, and one easily checks that
	\[l_\omega=\sum_{i=0}^k\sum_{(\omega_1',\omega_2')\in P_i}(-2)^{|\omega_1'|+|\omega_2'|},\]
	since $z_{\sigma_1}y_{\omega_1'}\cdot z_{\sigma_2}y_{\omega_2'}=z_{\sigma_1\cup\sigma_2}y_{\omega}$ for each $(\omega_1',\omega_2')\in P_i$ and the coefficient of $z_{\sigma_j}y_{\omega_j'}$ in $h_{\sigma_j,\omega_j}$ is $(-2)^{|\omega_j'|}$, $j=1,2$. Since $|\omega_1'|+|\omega_2'|=|\omega|+|\omega_1'\cap\omega_2'|$, a direct calculation gives the desired equation 
	\[l_\omega=(1-1)^k(-2)^{k+|\omega|}=\begin{cases}
		(-2)^{|\omega|},&k=0,\\
		0,&k>0,
	\end{cases}\] 
	and the proof is complete.
\end{proof}

\begin{proof}[Proof of Lemma \ref{lem:coboundary}]
	If $\omega_1\cap\omega_2=\emptyset$, there is nothing to prove. After relabeling the vertices, we may assume that $\omega_1\cap\omega_2=[s]$. Also, we use the convention that every face of $K$ is oriented by using the increasing order of its vertices.
	 For $i=1,2$, $j=1,2,\dots,s$, set 
	\begin{gather*}
	V_{i,j}=\{\sigma\in K_{\omega_i}\mid j\in\sigma,\,\sigma\cap\{1,\dots,j-1\}=\emptyset\},\\ 
	c_{i,j}=\sum_{\sigma\in V_{i,j}}l_{i,\sigma}\sigma^*,\quad
	b_{i,j}=c_i'+\sum_{k=j}^sc_{i,k},\\
	\begin{split}
	\alpha_j&=\psi_{\omega_1}(b_{1,j})\psi_{\omega_2}(b_{2,j})-\psi_{\omega_1}(b_{1,j+1})\psi_{\omega_2}(b_{2,j+1})\ \ (\text{set}\ b_{i,s+1}=c_i')\\
	&=\psi_{\omega_1}(c_{1,j})\psi_{\omega_2}(b_{2,j})+\psi_{\omega_1}(b_{1,j})\psi_{\omega_2}(c_{2,j})-\psi_{\omega_1}(c_{1,j})\psi_{\omega_2}(c_{2,j}).
	\end{split}
	\end{gather*}
Since $b_{i,1}=c_i$ and $b_{i,s+1}=c_i'$, $\alpha=\sum_{j=1}^s\alpha_j$. Hence, it suffices to prove that $2\alpha_j$ is a coboundary in $\QQ(K,\Lambda)$ for all $j$.

With the notation of \eqref{eq:notation}, for $i\in[m]$ let 
\[y_i=y_{0,i},\ \ z_i=x_{1,i}-y_{1,i},\ \ z_i'=y_{1,i},\ \ w_i=x_{2,i}+y_{2,i}.\]
Then by \eqref{eq:relations}, we have 
\begin{gather}
	y_i^2=y_i,\ \ y_iz_i=-z_i',\ \ z_iy_i=z_i+z_i',\ \ z_i^2=-w_i,\label{eq:produt}\\	
	dy_i=z_i,\ \ dz_i'=w_i.\label{eq:differential of y,z}
\end{gather}
From the definition of $\alpha_j$ and \eqref{eq:produt}, $\alpha_j$ can be written as 
\begin{equation}\label{eq:alpha_j}
	\alpha_j=w_jq_j+z_jr_j+z_j'r_j',\ \ 1\leq j\leq s,
\end{equation}
where $q_j$, $r_j$, $r_j'$ are defined as follows. Let $k_{i,\sigma}=(-1)^{|\sigma|}l_{i,\sigma}$, and let
\[W_{i,j}=\{\sigma\in K_{\omega_i}\mid\sigma\cap[s]=\emptyset\}\cup\bigcup_{k=j}^s V_{i,k}.\] 
Then
\begin{equation}\label{eq:q_j}
q_j=(\sum_{\sigma\in V_{1,j}}l_{1,\sigma} h_{\sigma\setminus\{j\},\omega_1\setminus\{j\}})(\sum_{\sigma\in V_{2,j}}k_{2,\sigma}h_{\sigma\setminus\{j\},\omega_2\setminus\{j\}}),
\end{equation}
\begin{equation}\label{eq:r_j}
\begin{split}
r_j&=(\sum_{\sigma\in V_{1,j}}k_{1,\sigma}h_{\sigma\setminus\{j\},\omega_1\setminus\{j\}})(\sum_{\sigma\in W_{2,j+1}}k_{2,\sigma}\sum_{j\in\omega'\subset\omega_2\setminus\sigma}-(-2)^{|\omega'|-1}z_\sigma y_{\omega'\setminus\{j\}})\\
&+(\sum_{\sigma\in W_{1,j+1}}l_{1,\sigma} \sum_{j\in\omega'\subset\omega_1\setminus\sigma}(-2)^{|\omega'|-1}z_\sigma y_{\omega'\setminus\{j\}})(\sum_{\sigma\in V_{2,j}}k_{2,\sigma}h_{\sigma\setminus\{j\},\omega_2\setminus\{j\}}),
\end{split}
\end{equation}
\begin{equation}\label{eq:r_j'}
\begin{split}
	r_j'&=(\sum_{\sigma\in V_{1,j}}k_{1,\sigma} h_{\sigma\setminus\{j\},\omega_1\setminus\{j\}})(\sum_{\sigma\in W_{2,j+1}}k_{2,\sigma}\sum_{j\in\omega'\subset\omega_2\setminus\sigma}(-2)^{|\omega'|}z_\sigma y_{\omega'\setminus\{j\}})\\
&-(\sum_{\sigma\in W_{1,j+1}}l_{1,\sigma}\sum_{j\in\omega'\subset\omega_1\setminus\sigma}(-2)^{|\omega'|}z_\sigma y_{\omega'\setminus\{j\}})(\sum_{\sigma\in V_{2,j}}k_{2,\sigma} h_{\sigma\setminus\{j\},\omega_2\setminus\{j\}}).
\end{split}
\end{equation}
These formulas come from \eqref{eq:map phi_omega}, \eqref{eq:produt} and the definitions of $\alpha_j$ and $h_{\sigma,\omega}$. \eqref{eq:r_j} also uses the equality
\[
	\sum_{j\in\omega'\subset\omega_i\setminus\sigma}(-2)^{|\omega'|}z_\sigma y_{\omega'}=y_j\sum_{j\not\in\omega'\subset\omega_i\setminus\sigma}(-2)^{|\omega'|+1}z_\sigma y_{\omega'},\ \ j\not\in\sigma.
\]
From \eqref{eq:r_j} and \eqref{eq:r_j'}, we see that $r_j'=2r_j$.

The proof of $2\alpha_j\in d\QQ(K,\Lambda)$ is by induction on $j$, setting $\alpha_0=0$ as the base case. Assume that the result has been proved for all $\alpha_k$ with $k<j$. By Lemma \ref{lem:product of h}, $\psi_{\omega_1}(c_1')\psi_{\omega_2}(c_2')=\psi_{\omega_1+\omega_2}(c_1'*c_2')$.
Hence, $d(\psi_{\omega_1}(c_1')\psi_{\omega_2}(c_2'))=0$ since $c_1'*c_2'$ is a cocycle in $\w C^*(K_{\omega_i+\omega_2})$. It follows that $d\alpha=0$, and then
 $d(\sum_{k=j}^s\alpha_k)=0$ since $d(\sum_{k=1}^{j-1}\alpha_k)=0$ by the induction hypothesis (note that $\QQ_K$ is torsion free, so $2d\alpha_k=0$ if and only if $d\alpha_k=0$). By \eqref{eq:differential of y,z} and \eqref{eq:alpha_j}, 
\[
	d\alpha_k=w_k(dq_k+r_k')-z_kdr_k-z_k'dr_k',\ \ 1\leq k\leq s.
\]
This implies that 
\[dq_j+r_j'\in E_j:=\{f\in\QQ_K\mid w_jf=0\},\]  
since $r_j,r_j'$ are polynomials in variables $y_i,z_i$ with $i\neq j$, and for each $k>j$, $\alpha_k$ is a polynomial in variables $y_i$ with $i\in[m]$ and $z_i$ with $i>j$, so that the $j$th degrees (in the $\Zz^m$-grading) of $z_jdr_j$, $z_j'dr_j'$ and $d\alpha_k$ for $k>j$ are at most $1$. Using \eqref{eq:relations} and the definition of $\QQ_K$, one can show that $E_j$ is spanned by basis elements $e=e_1\otimes\cdots\otimes e_m\in\QQ_K$ such that $\supp e\not\in\st_K\{j\}$ and that 
\[E_j=\{f\in\QQ_K\mid z_jf=0\}\subset\{f\in\QQ_K\mid z_j'f=0\}.\]

Set $\eta_j=(2z_j'+z_j)q_j$. We claim that $\eta_j\in\QQ(K,\Lambda)$. Indeed, given $\sigma,\tau\in K$ with $j\in\sigma\subset\omega_1$ and $j\in\tau\subset\omega_2$, let $p_{\sigma,\tau}=h_{\sigma\setminus\{j\},\omega_1\setminus\{j\}}h_{\tau\setminus\{j\},\omega_2\setminus\{j\}}$. Then 
\[z_j^2p_{\sigma,\tau}=\pm h_{\sigma,\omega_1}h_{\tau,\omega_2}.\] 
Since $h_{\sigma,\omega_1},h_{\tau,\omega_2}\in\QQ(K,\Lambda)$ by Lemma \ref{lem:invariants} and $z_i^2=-(x_{2,i}+y_{2,i})\in\QQ(K,\Lambda)$ for any $i\in[m]$, $z_j^2g(p_{\sigma,\tau})=z_j^2p_{\sigma,\tau}$ for any $g\in\ker\Lambda$, and so $gp_{\sigma,\tau}-p_{\sigma,\tau}\in E_j$.
Using \eqref{eq:q_j}, this implies that $gq_j-q_j\in E_j$ for any $g\in\ker\Lambda$. Furthermore, since $2z_j'+z_j=x_{1,j}+y_{1,j}\in\QQ(K,\Lambda)$, it follows that 
\[g\eta_j-\eta_j=(2z_j'+z_j)(gq_j-q_j)=0 \text{ for any } g\in\ker\Lambda.\] 
Hence our claim is true.
The previous arguments show that
\[d\eta_j=2w_jq_j-z_jdq_j-2z_j'dq_j=2w_jq_j+2z_jr_j+2z_j'r_j'=2\alpha_j,\] 
finishing the induction step.
\end{proof}

\begin{proof}[Proof of Theorem \ref{thm:ring isomorphism}]
	Following the notation of Lemma \ref{lem:coboundary} and its proof, let $\eta=\sum_{j=1}^s\eta_j\in\QQ(K,\Lambda)$, and let $\bar\eta\in\QQ(K,\Lambda)\otimes\Zz/2$ be the mod $2$ reduction of $\eta$.  Note that $h_{\sigma,\omega}\equiv z_{\sigma}$ mod $2$. It follows that
	\[\bar\eta=\sum_{j=1}^s\big(z_j(\sum_{\sigma\in V_{1,j}}l_{1,\sigma}z_{\sigma\setminus\{j\}})(\sum_{\sigma\in V_{2,j}}l_{2,\sigma} z_{\sigma\setminus\{j\}})\big).\]
	By Theorem \ref{thm:DGA}, there are ring isomorphisms
	\begin{gather*}
	\iota:H^*(\QQ(K,\Lambda))\to H^*(M(K,\Lambda)),\\
	\iota_{\Zz/2}:H^*(\QQ(K,\Lambda)\otimes\Zz/2)\to H^*(M(K,\Lambda);\Zz/2).
     \end{gather*}
	As we saw in the proof of Theorem \ref{thm:mod 2 image}, $\iota_{\Zz/2}(z_\sigma)=x_\sigma$, using \eqref{diag:Phi and Psi}. Thus $\iota_{\Zz/2}(\bar\eta)=E(c_1,c_2)$. Since $d\eta=2\alpha$, we have $\beta_{\Zz}([\bar\eta])=[\alpha]$. Combining these together, we see that
	\[\beta_{\Zz}(E(c_1,c_2))=\iota\left(\Psi_\Lambda([c_1])\Psi_\Lambda([c_2])-\Psi_\Lambda([c_1'*c_2'])\right).\]
	
	From the proof of Proposition \ref{prop:group isomorphism}, we know that the map $\bar\kappa:\A_2/J\to S/I$ is the composition $\bar\beta_{\Zz}^{-1}\circ\iota\circ\Psi_\Lambda$.
	Hence, the theorem will follow once we prove that the element $E(c_1,c_2)\in S$ 
	does not depend on the choice of the representatives of $[c_1]$ and $[c_2]$.
	Suppose that $c_i$ and $\tilde c_i$, $i=1,2$, represent the same class in $H^*(\A(K,\Lambda))$. Since the differential of $\A(K,\Lambda)$ is $d=2\delta$, there is a cochain $b_i$ such that $\tilde c_i=c_i+2\delta b_i$.
	Moreover, $E(c_1,c_2)$ depends only on the coefficients of $c_1,c_2$ modulo $2$. Thus, $E(c_1,c_2)=E(\tilde c_1,\tilde c_2)$ and the proof is complete. 
\end{proof}

\bibliography{M-A}
\bibliographystyle{amsplain}

\end{document}